\documentclass[12pt]{article}
\usepackage{amsmath,amsfonts,amssymb}
\newtheorem{definition}{Definition}[section]

\newtheorem{proposition}{Proposition}[section]
\newtheorem{property}{Property}[section]
\newtheorem{theorem}{Theorem}[section]

\def\cA{{\cal A}}

\def\cP{{\cal P}}                    
\def\cS{{\cal S}}          \def\cT{{\cal T}}          \def\cU{{\cal U}}
          \def\cW{{\cal W}}

\def\CC{{\mathbb C}}
\def\HH{{\mathbb H}}
\def\II{{\mathbb I}}
\def\ZZ{{\mathbb Z}}

\def\tr{\mathop{\rm Tr}\nolimits}
\def\qdet{\text{q-}\!\det}
\def\eps{\varepsilon}

\newcommand{\sfrac}[2]{{\textstyle{\frac{#1}{#2}}}}
\newcommand{\half}{{\scriptstyle{\frac{1}{2}}}}
\newcommand{\three}{{\scriptstyle{\frac{3}{2}}}}
\newcommand{\car}[2]{\genfrac{[}{]}{0pt}{}{#1}{#2}}
\newcommand{\finproof}{{\hfill \rule{5pt}{5pt}\\}}

\numberwithin{equation}{section}

\begin{document}
\baselineskip=18pt
\pagestyle{empty}
\pagenumbering{arabic}
\setcounter{page}{0}
\parindent=0pt

\markright{DRAFT\dotfill\today\dotfill}

\null
\vfill

\begin{center}

{\Large \textsf{Sugawara and vertex operator constructions\\[2.1ex]
for deformed Virasoro algebras}}

\vspace{10mm}

{\large D. Arnaudon$^a$,
  J. Avan$^{b}$\footnote{avan@ptm.u-cergy.fr, 
  frappat@lapp.in2p3.fr, ragoucy@lapp.in2p3.fr, shiraish@ms.u-tokyo.ac.jp\label{foot:1}},  
  L. Frappat$^{ac\ref{foot:1}}$,
  {E}. Ragoucy$^{a\ref{foot:1}}$, J. Shiraishi$^{d\ref{foot:1}}$} 

\vspace{10mm}

\emph{$^a$ Laboratoire d'Annecy-le-Vieux de Physique Th{\'e}orique}

\emph{LAPTH, CNRS, UMR 5108, Universit{\'e} de Savoie}

\emph{B.P. 110, F-74941 Annecy-le-Vieux Cedex, France}

\vspace{7mm}

\emph{$^b$ Laboratoire de Physique Th{\'e}orique et Mod{\'e}lisation}

\emph{Universit{\'e} de Cergy, 5 mail Gay-Lussac, Neuville-sur-Oise}

\emph{F-95031 Cergy-Pontoise Cedex}

\vspace{7mm}

\emph{$^c$ Member of Institut Universitaire de France}

\vspace{7mm}

\emph{$^d$ Graduate School of Mathematical Science, }

\emph{University of Tokyo, Komaba, Meguro-ku, Tokyo,
153-8914, Japan}
\end{center}

\vfill

\begin{abstract}
From the  defining exchange relations of the ${\cA}_{q,p}(\widehat{gl}_{N})$
 elliptic quantum algebra,
we construct subalgebras which can be characterized
as $q$-deformed $W_N$ algebras. The consistency conditions relating
the parameters $p,q,N$ and the central charge $c$ are shown to be related to the
singularity structure of the functional coefficients
defining the exchange relations of 
specific vertex operators  representations 
of ${\cA}_{q,p}(\widehat{gl}_{N})$ available when $N=2$.
\end{abstract}
\vfill
\centerline{\textit{Dedicated to our friend Daniel Arnaudon}}
\vfill
MSC: 81R10, 17B37, 17B69 ---
PACS: 02.20.Uw, 02.20.-c, 02.30.Ik
\vfill

\rightline{math.QA/0601250}
\rightline{LAPTH-1135/06}
\rightline{January 2006}
\newpage
\pagestyle{plain}

\section{Introduction}

The notion of $q$-deformed Virasoro, and its natural extensions the
$q$-deformed $W_n$ algebras, covers a number of algebraic structures
characterized by quadratic exchange relations such that one can define
a semi-classical limit as a Poisson structure with a very specific
form.  This latter exhibits in turn a deformation parameter $q$ such
that the suitably defined limit $ q \to 1$ leads to the usual
classical Virasoro and $W_n$ Poisson algebras.  The original Poisson
structure was proposed in \cite{FR} as originating on an extended
center of the affine quantum algebra $U_q(\widehat{sl}_{2})_{c}$ for
the critical value $c=-2$ of the central charge.  Quantization of this
structure has proved to be an intringuing task \cite{SKAO,FF,AKOS95}.  
Indeed, it has been found
deep relationship among the Virasoro and $W_n$ algebras,
algebras of screening currents, the quantum Miura transformation, and
the Macdonald polynomials or the
quantum $N$-body trigonometric Ruijsenaars--Schneider model.

Another connection between $q$-deformed Virasoro and $W_N$ algebras, and the
vertex-type elliptic quantum algebras ${\cA}_{q,p}(\widehat{sl}_{2})$ of
\cite{JKOS}, was established in a series of papers
\cite{ENSLAPP649,ENSLAPP670,LAPTH690}. It was shown that the formal
defining exchange relations of the vertex elliptic algebras, parametrized
by elliptic $R$ matrices of Baxter--Belavin \cite{Bela}, allowed to define
certain operators, bilinear in the quantum Lax matrix, realizing 
quadratic exchange relations of $q$-$W_N$ type (i.e. with suitable
semi-classical limits), whenever particular commensurability relations
existed between the elliptic module $p$, the deformation parameter $q$ and
the central charge $c$. Supplementary commensurability relations could then
lead to commuting operators, allowing for a definition of consistent
Poisson structures reproducing the classical $q$-$W_N$ algebra. 

At that point, it may be useful to make a comparison with the constructions
done in the undeformed case. The first approach (based on MacDonald
polynomials and deformed bosons) can be seen as a deformation of the Miura
construction for Virasoro and $W_{N}$ algebras \cite{luky},
which realizes these algebras in terms of free fields. The second approach
(based on quadratics of ${\cA}_{q,p}(\widehat{sl}_{2})$ generators) can be
viewed as a deformation of the Sugawara construction
\cite{sug,KZ1984,Godol85}, which realizes the Virasoro and $W_{N}$ algebras
in terms of affine Lie algebras. In the undeformed case, these two
constructions lead to the same algebras, and their connection is done when
considering Sugawara construction for an affine Lie algebra, itself
realized in terms of vertex operators. When considering the deformations of
these approaches, one gets \textit{two} types of algebras, with different
structure constants, but same classical limit. Moreover, these
constructions have remained until now at the level of formal manipulations
since explicit realizations of the elliptic algebra were lacking at the
time. The main objective of this paper is to present a first contact
point between these two approaches at the quantum level.

Indeed, the recent construction of explicit vertex operator (V.O.) representations
of the elliptic vertex algebra (for particular values of $c$ and relations
between $p$ and $q$), and the identification of specific fused operators in
this vertex algebra \cite{Junichi04}, lead us to reconsider and
generalize our former construction of $q$-$W_N$ generators from this point
of view. New sets of commensurability relations have now been identified,
both for ${\cA}_{q,p}(\widehat{sl}_{2})$ and
${\cA}_{q,p}(\widehat{sl}_{N})$. In addition it is now possible to express
our formal bilinear operators in terms of this explicit vertex algebra
construction, at least for $sl(2)$, and therefore to examine the meaning of
the commensurability relations also in the framework of this V.O.
representation. It turns out as we shall now see, that these relations are
interpreted as characterizing locii for simultaneous singularities in the
exchange relations of the two types of V.O. The bilinear operators
generating our $q$-$W_N$ algebras are then represented as products of
residues of both V.O. at these singular points.

The plan of the paper is as follows. After a reminder of the definition and
some useful properties of ${\cA}_{q,p}(\widehat{gl}_{2})$ in section
\ref{sect:elp}, we define, in section \ref{sect:exch}, operators labelled
by two integers $\ell$ and $\ell'$ realizing an
${\cA}_{q,p}(\widehat{gl}_{2})$ subalgebra when the parameters $p,q$ and
the central charge $c$ obey certain consistency conditions. The
semi-classical limit yields the classical deformed Virasoro algebra
introduced by \cite{FR}. In section \ref{sect:VO}, we recall the level one
vertex operator representation of ${\cA}_{q,p}(\widehat{gl}_{2})$
introduced by \cite{FIJKMY95,Junichi04} and we connect these vertex
operators with the above subalgebra. The generalization to
${\cA}_{q,p}(\widehat{gl}_{N})$ is considered in section \ref{sect:glN}.

\section{The elliptic algebra ${\cA}_{q,p}(\widehat{gl}_{2})$\label{sect:elp}}

The quantum affine elliptic algebra ${\cA}_{q,p}(\widehat{gl}_{2})$ is
described in the RLL formalism \cite{FIJKMY}. Its R-matrix
is constructed from the Boltzmann weights of the eight vertex model:
\begin{equation}
  \label{eq:RAqp}
  \cW(z,q,p) = \rho(z) \left(
  \begin{array}{cccc}
    a(z) & 0 & 0 & d(z) \\
    0 & b(z) & c(z) & 0 \\
    0 & c(z) & b(z) & 0 \\
    d(z) & 0 & 0 & a(z) \\
  \end{array}
  \right)
\end{equation}
with
\begin{equation}
  \label{eq:abcdAqp}
  \begin{array}{lcl}
    a(z) = z^{-1} \; \dfrac{\Theta_{p^2}(q^2z^2) \;
    \Theta_{p^2}(pq^2)} {\Theta_{p^2}(pq^2z^2) \; \Theta_{p^2}(q^2)}
    &\qquad& d(z) = -\dfrac{p^{\half}}{q z^{2}} \;
    \dfrac{\Theta_{p^2}(z^2) \; \Theta_{p^2}(q^2z^2)}
    {\Theta_{p^2}(pz^2) \; \Theta_{p^2}(pq^2z^2)} \\[5mm] 
    b(z) = qz^{-1} \; \dfrac{\Theta_{p^2}(z^2) \;
    \Theta_{p^2}(pq^2)} {\Theta_{p^2}(pz^2) \; \Theta_{p^2}(q^2)}
    &\qquad& c(z) = 1 \\
  \end{array}
\end{equation}
with $\Theta_{p}(z) = (z;p)_{\infty}\, (pz^{-1};p)_{\infty} \,
(p;p)_{\infty}$ and $(z;a_1,\dots,a_m)_\infty = \prod_{n_i \ge 0} (1-z
a_1^{n_1} \dots a_m^{n_m})$. \\
The normalization factor $\rho(z)$ is chosen as follows:
\begin{equation}
  \label{eq:normAqp}
  \rho(z) = \frac{(p^2;p^2)_\infty} {(p;p)_\infty^2} \;
  \frac{\Theta_{p^2}(pz^2)\Theta_{p^2}(q^2)}
  {\Theta_{p^2}(q^2z^2)} \;
  \frac{\xi(z^2;p,q^4)}{\xi(z^{-2};p,q^4)}
\end{equation}
where we have introduced
\begin{equation}
  \label{eq:xi}
  \xi(z;p,q^4) = \frac{(q^2z;p,q^4)_\infty \;
  (pq^2z;p,q^4)_\infty} {(q^4z;p,q^4)_\infty \;
  (pz;p,q^4)_\infty}\,.
\end{equation}

\begin{property}
  The matrix (\ref{eq:RAqp}) has the following properties:
  \begin{align}
    & \text{YBE:} && \cW_{12}(z_1) \,
    \cW_{13}(z_2) \, \cW_{23}(z_2/z_1) =
    \cW_{23}(z_2/z_1) \, \cW_{13}(z_2) \, \cW_{12}(z_1)
    & \\[6pt] 
    & \text{unitarity:} && \cW_{12}(z) \, \cW_{21}(z^{-1}) = 1 &
    \\[6pt] 
    & \text{crossing symmetry:} && \cW_{21}^{t_1}(\frac{1}{z}) = (\sigma_{x}
    \otimes 1) \cW_{12}(-\frac{z}{q}) (\sigma_{x} \otimes 1) = (1 \otimes
    \sigma_{x}) \cW_{12}(-\frac{z}{q}) (1 \otimes \sigma_{x}) & \\[6pt]
    & \text{antisymmetry:} && \cW_{12}(-z) = - (\sigma_{z} \otimes 1)
    \cW_{12}(z) (\sigma_{z} \otimes 1) = - (1 \otimes \sigma_{z})
    \cW_{12}(z) (1 \otimes \sigma_{z}) &
  \end{align}
  where $\sigma_{x},\sigma_{y},\sigma_{z}$ are the $2 \!\times\! 2$ Pauli
  matrices and $t_i$ denotes the transposition in space $i$.
\end{property}

\bigskip

The definition of the quantum affine elliptic algebra
${\cA}_{q,p}(\widehat{gl}_{2})$ requires the use of a slightly modified
R-matrix $R_{12}(z)$, which differs from (\ref{eq:RAqp}) by
a suitable normalization factor:
\begin{equation}
  R_{12}(z,q,p)\equiv R_{12}(z)
  = \tau(q^{\half}z^{-1}) \cW_{12}(z,q,q)\,.
\end{equation}
The factor $\tau(z)$ is given by
\begin{equation}
  \tau(z) = z^{-1} \,
  \frac{\Theta_{q^4}(qz^2)}{\Theta_{q^4}(qz^{-2})}
\end{equation}
which is $q^2$-periodic and satisfies $\tau(z) \tau(z^{-1}) = 1$ and
$\tau(q^{\half}z^{-1}) = -\tau(q^{\half}z)$. \\
The R-matrix $R_{12}(z)$ is no longer unitary but verifies
\begin{equation}
  R_{12}(z) \, R_{21}(z^{-1}) = \widetilde{\tau}(z)
  \label{eq:unitRtilde}
\end{equation}
where the $\widetilde{\tau}$ function is given by
\begin{equation}
    \widetilde{\tau}(z) = \tau(q^{\half}z^{-1}) \, \tau(q^{\half}z) = -
    \tau(q^{\half}z)^2 = - \left[ q^{-\half}z \;
    \frac{\Theta_{q^4}(q^2z^2)}{\Theta_{q^4}(z^2)} \right]^2
    \,.\label{eq:tautilde}
\end{equation}
The R-matrix $R_{12}(z)$ obeys a quasi-periodicity property
\begin{equation}
  \label{eq:quasiper}
  R_{12}(-p^{\half}z) = (\sigma_{x} \otimes 1) \left(
  R_{21}(z^{-1}) \right)^{-1} (\sigma_{x} \otimes 1) =
  \widetilde{\tau}(z)^{-1} \, (\sigma_{x} \otimes 1) \,
  R_{12}(z) \, (\sigma_{x} \otimes 1)\,,
\end{equation}
so that a recursive use of this formula leads to
\begin{equation}
  \label{eq:relF}
  R_{12}\left((-p^{\half})^{-\ell}z \right) = F(\ell,z)
  \, (\sigma_{x}^\ell \otimes 1) \, R_{12}(z) \,
  (\sigma_{x}^\ell \otimes 1)
\end{equation}
where
\begin{equation}
  \label{eq:Fll}
  F(\ell,z) = \left\{
  \begin{array}{cl}
    \displaystyle \prod_{s=1}^\ell
    \widetilde{\tau}\big((-p^\half)^{-s}z\big) & \mbox{ for
	$\ell\geq 0$}
    \\
    \\ 
    \displaystyle \prod_{s=1}^{|\ell|}
    \widetilde{\tau}\big((-p^\half)^{|\ell|-s}z\big)^{-1} & \mbox{ for
    $\ell \leq 0$}
  \end{array}
  \right.
\end{equation}
The function $F(\ell,z)$ is $q^2$-periodic (due to the
$q^2$-periodicity of the $\tau$ function) and satisfies the following
relations ($\ell, n \in \ZZ$)
\begin{eqnarray}
  \label{eq:pptF1}
  F(\ell,z) F(-\ell,z^{-1}) =
  \frac{\widetilde{\tau}\big((-p^{\half})^{-\ell}z\big)}
  {\widetilde{\tau}(z)} \\
  \label{eq:pptF2}
  F(\ell,(-p^{\half})^{-n}z) = \frac{F(\ell+n,z)}{F(n,z)} .
\end{eqnarray}
In particular, one has 
\begin{equation}
  \label{eq:inver}
  F(-\ell,z) = \frac{1}{F(\ell,(-p^{\half})^{\ell}z)}.
\end{equation}
The crossing symmetry and the unitarity properties of $\cW_{12}$ then allow
one to exchange inversion and transposition for the matrix
$R_{12}$:
\begin{equation}
  \label{eq:crossunit}
  \Big( R_{12}(z)^{t_2} \Big)^{-1} = \Big(
  R_{12}(q^2z)^{-1} \Big)^{t_2}\,.
\end{equation}

\medskip

The quantum affine elliptic algebra ${\cA}_{q,p}(\widehat{gl}_{2})$ is
defined as a formal algebra of operators
\begin{equation}
  L(z) = \sum_{n\in\ZZ} L_{n} \, z^{-n} = \left(
  \begin{array}{cc}
    L_{++}(z) & L_{+-}(z) \\
    L_{-+}(z) & L_{--}(z) \\
  \end{array}
  \right) = \sum_{i,j = \pm} L_{ij}(z) E_{ij}
\end{equation}
where the functions $L_{++}$ and $L_{--}$ are even while $L_{+-}$ and
$L_{-+}$ are odd in the variable $z$, and obey the following relations:
\begin{equation}
  \label{eq:RLLAqp}
  R_{12}(z_1/z_2) \, L_{1}(z_1) \, L_{2}(z_2)
  = L_{2}(z_2) \, L_{1}(z_1) \,
  R^{*}_{12}(z_1/z_2)
\end{equation}
with $L_1(z) = L(z) \otimes 1$, $L_2(z) = 1 \otimes L(z)$
and $R^{*}_{12}(z) \equiv
R_{12}(z,q,p^*=pq^{-2c})$. \\
The quantum determinant of $L(z)$ given by \cite{FIJKMY95}
\begin{equation}
  \qdet L(z) = L_{++}(q^{-1}z) L_{--}(z) -
  L_{-+}(q^{-1}z) L_{+-}(z)
\end{equation}
is in the center of ${\cal A}_{q,p}(\widehat{gl}_{2})$. It can be factored
out, and set to the value $q^{\frac c2}$ ($c$ being the central charge) so
as to get
\begin{equation}
    {\cA}_{q,p}(\widehat{sl}_{2}) = \cA_{q,p}(\widehat{gl}_{2})/
    \big\langle \qdet L(z) - q^{\frac{c}{2}} \big\rangle\,.
\end{equation}

\section{Exchange algebras\label{sect:exch}}

\subsection{Deformed Virasoro algebras\label{sect:DVA}}

\begin{proposition}\label{prop:LT}
  The operators $T_{\ell\ell'}(z)$ defined by ($\ell,\ell'
  \in \ZZ$)
  \begin{equation}
    \label{eq:operT}
    T_{\ell\ell'}(z) \equiv \tr \Big( L(\gamma z)^{-1} \;
    \sigma_{x}^{\ell} \; L(z) \; \sigma_{x}^{\ell'} \Big) = \tr \Big(
    \sigma_{x}^{\ell} \; {(L(\gamma z)^{-1})}^{t} \; \sigma_{x}^{\ell'}
    \; {L(z)}^{t} \Big)
  \end{equation}
   have exchange relations
  with the generators $L(z)$ of ${\cA}_{q,p}(\widehat{gl}_{2})$:
  \begin{equation}
    L(z_2) \; T_{\ell\ell'}(z_1) = f_{\ell\ell'}(z_1/z_2)
    \; T_{\ell\ell'}(z_1) \; L(z_2)
  \end{equation}
  if the conditions 
 \begin{equation}
   \gamma = \gamma_{\ell\ell'} \equiv (-p^{\half})^{-\ell} =
   (-p^{*\half})^{-\ell'}\, q^{2}
\label{gammaT}
\end{equation}
are fulfilled. \\
The exchange function $f_{\ell\ell'}(z)$ is
  given by
  \begin{equation}
    \label{eq:functstrucTL}
    f_{\ell\ell'}(z) = \frac{F^*(\ell',z)}{F(\ell,z)}
  \end{equation}
  where $F(\ell,z)$ is given by (\ref{eq:Fll}) and $F^{*}(\ell,z)$
  is obtained from $F(\ell,z)$ through the shift
  $p \to p^{*} = pq^{-2c}$.
\end{proposition}

\textbf{Proof:} One has
\begin{eqnarray}
 && L_{2}(z_2) \, {T_{1}}_{\ell\ell'}(z_1) = \tr_1 \Big(
  L_{2}(z_2) \; \sigma_{1x}^{\ell} \; {(L_{1}(\gamma
  z_1)^{-1})}^{t_{1}} \; \sigma_{1x}^{\ell'} \;
  {L_{1}(z_1)}^{t_{1}} \Big) \nonumber \\
  &&= \tr_1 \Big( \sigma_{1x}^{\ell} \;
  {(R_{12}^{-1}(\gamma z_1/z_2))}^{t_{1}} \;
  {(L_{1}(\gamma z_1)^{-1})}^{t_{1}} \; L_{2}(z_2) \;
  {({(R_{12}^{*-1}(\gamma z_1/z_2))}^{t_{1}})}^{-1} \;
  \sigma_{1x}^{\ell'} \; {L_{1}(z_1)}^{t_{1}} \Big) \nonumber \\
  &&= \tr_1 \Big( \sigma_{1x}^{\ell} \;
  {(R_{12}^{-1}(\gamma z_1/z_2))}^{t_{1}} \;
  {(L_{1}(\gamma z_1)^{-1})}^{t_{1}} \; L_{2}(z_2) \;
  {(R_{12}^{*}(q^{-2} \gamma z_1/z_2))}^{t_{1}} \;
  \sigma_{1x}^{\ell'} \; {L_{1}(z_1)}^{t_{1}} \Big) \nonumber 
\end{eqnarray}
where the crossing-unitarity property (\ref{eq:crossunit}) has been used in
the last equality. \\
Imposing the relation $q^{-2} \gamma = (-p^{*\half})^{-\ell'}$ and using
the quasi-periodicity property (\ref{eq:relF}) of $R$ allows
one to further exchange $L_{2}(z_2)$ and $L_{1}(z_1)^{t_{1}}$:
\begin{eqnarray}
  L_{2}(z_2) \, {T_{1}}_{\ell\ell'}(z_1) &=&
  F^{*}(\ell',z_1/z_2) \, \tr_1 \Big( \sigma_{1x}^{\ell} \;
  {(R_{12}^{-1}(\gamma z_1/z_2))}^{t_{1}} \;
  {(L_{1}(\gamma z_1)^{-1})}^{t_{1}} \; \sigma_{1x}^{\ell'} \;
  \nonumber \\
  && \hspace*{20ex} L_{2}(z_2) \;
  {(R_{12}^{*}(z_1/z_2))}^{t_{1}} \;
  {L_{1}(z_1)}^{t_{1}} \Big) \nonumber \\
  &=& F^{*}(\ell',z_1/z_2) \, \tr_1 \Big( \sigma_{1x}^{\ell} \;
  {(R_{12}^{-1}(\gamma z_1/z_2))}^{t_{1}} \;
  {(L_{1}(\gamma z_1)^{-1})}^{t_{1}} \; \sigma_{1x}^{\ell'} \; \nonumber
  \\
  && \hspace*{20ex} {L_{1}(z_1)}^{t_{1}} \;
  {(R_{12}(z_1/z_2))}^{t_{1}} \; L_{2}(z_2) \Big)
\end{eqnarray}
Imposing now the relation $\gamma = (-p^{\half})^{-\ell}$ and using once
again the relation (\ref{eq:relF}), one gets
\begin{eqnarray}
  L_{2}(z_2) \, {T_{1}}_{\ell\ell'}(z_1) &=&
  F^*(\ell',z_1/z_2) \, F(\ell,z_1/z_2)^{-1} \, \tr_1 \Big(
  {(R_{12}^{-1}(z_1/z_2))}^{t_{1}} \; \sigma_{1x}^{\ell}
  \; {(L_{1}(\gamma z_1)^{-1})}^{t_{1}} \; \sigma_{1x}^{\ell'} \;
  \nonumber \\
  && \hspace*{50mm} {L_{1}(z_1)}^{t_{1}} \;
  {(R_{12}(z_1/z_2))}^{t_{1}} \; L_{2}(z_2) \Big)
\end{eqnarray}
Using the fact that under a trace over the space 1 one has ${\rm Tr}_1
\Big( R_{12} Q_1 {R'}_{12} \Big) = {\rm Tr}_1 \Big( Q_1 {R'_{12}}^{t_2}
{R_{12}}^{t_2} \Big)^{t_2}$, one gets
\begin{eqnarray}
  L_{2}(z_2) \, {T_{1}}_{\ell\ell'}(z_1) &=&
  F^*(\ell',z_1/z_2) \, F(\ell,z_1/z_2)^{-1} \, \tr_1 \Big(
  \sigma_{1x}^{\ell} \; {(L_{1}(\gamma z_1)^{-1})}^{t_{1}} \;
  \sigma_{1x}^{\ell'} \; {L_{1}(z_1)}^{t_{1}} \; \nonumber \\
  && \hspace*{30mm} {(R_{12}(z_1/z_2))}^{t_{1}t_{2}}
  \; {(R_{12}^{-1}(z_1/z_2))}^{t_{1}t_{2}} \;
  \Big)^{t_2} \; L_{2}(z_2)
\end{eqnarray}
which leaves a trivial dependence in space 2 under the trace in space 1 and
therefore leads to the exchange relation between
${T_{1}}_{\ell\ell'}(z_1)$ and $L_{2}(z_2)$:
\begin{equation}
  L_{2}(z_2) \; {T_{1}}_{\ell\ell'}(z_1) =
  F^*(\ell',z_1/z_2) \, F(\ell,z_1/z_2)^{-1} \;
  {T_{1}}_{\ell\ell'}(z_1) \; L_{2}(z_2)
\end{equation}
\finproof
Let us stress that although $T_{\ell\ell'}(z)$ seems to depend on
$\ell$ and $\ell'$ only modulo 2, it really depends on them because of 
the form of $\gamma$ imposed by conditions (\ref{gammaT}).

In the same way, one can establish:
\begin{proposition}\label{prop:LS}
  The operators $S_{\ell\ell'}(z)$ defined by ($\ell,\ell' \in \ZZ$)
  \begin{equation}
    \label{eq:operS}
    S_{\ell\ell'}(z) \equiv \tr \Big( \sigma_{x}^{\ell} \; L(z) \;
    \sigma_{x}^{\ell'} \; {(L(\gamma^{-1}q^{2}z)^{-1})} \Big) = \tr \Big(
    {L(z)}^{t} \; \sigma_{x}^{\ell} \; {(L(\gamma^{-1}q^{2}z)^{-1})}^{t}
    \; \sigma_{x}^{\ell'} \Big)
  \end{equation}
  have exchange relations
  with the generators $L(z)$ of ${\cA}_{q,p}(\widehat{gl}_{2})$:
  \begin{equation}
    L(z_2) \; S_{\ell\ell'}(z_1) = f_{\ell\ell'}(z_1/z_2)
    \; S_{\ell\ell'}(z_1) \; L(z_2)
  \end{equation}
  if the conditions 
  \begin{equation}
    \gamma = \gamma_{-\ell-\ell'} \equiv (-p^{\half})^{\ell} =
    (-p^{*\half})^{\ell'}\, q^{2}
    \label{gammaS}
\end{equation}
are fulfilled.
\end{proposition}
The proof follows along the same lines as for proposition \ref{prop:LT}.
\finproof

Let us remark that the above calculations rely only on the exchange
relation (\ref{eq:RLLAqp}). As such, these formal manipulations do not
allow us to determine potential `delta-type' terms. These
terms can be obtained through explicit realizations of the ${\cal A}_{q,p}(\widehat{gl}_{2})$
algebra. Unfortunately, they are still missing.

Note that the relations (\ref{gammaS}) are deduced from the relations 
(\ref{gammaT}) through the change $(\ell,\ell')\to(-\ell,-\ell')$. 
The above results suggest the following definition:
\begin{definition}
In the parameter space $(q,p,c)$, the surface $\cP_{\ell\ell'}$ is
defined by the relation 
\begin{equation}
(-p^{\half})^{-\ell} = (-p^{*\half})^{-\ell'}\,q^{2}\,,\mbox{ with }
p^{*}=p\,q^{-2c}
\end{equation}
The operators $T_{\ell\ell'}$ and $S_{-\ell,-\ell'}$ defined on this
surface are given by relations (\ref{eq:operT}) and (\ref{eq:operS}) with
$\gamma\equiv\gamma_{\ell\ell'}=(-p^{\half})^{-\ell}$.
\end{definition}
Since $T_{\ell\ell'}$ and $S_{\ell\ell'}$ have the same exchange relations
with the $L(z)$ operators, one can say that $T_{\ell\ell'}$ and
$S_{\ell\ell'}$ represent the same formal algebra on respectively the
surfaces $\cP_{\ell\ell'}$ and $\cP_{-\ell,-\ell'}$.

{From} Propositions \ref{prop:LT} and \ref{prop:LS}, we get immediately:
\begin{proposition}\label{prop:TT}
When $c$ is a rational number, several surfaces $\cP_{\ell\ell'}$ can 
be defined simultaneously (i.e. have non-trivial intersection). \\
  On the surface $\cP_{\ell\ell'}\cap\cP_{\lambda\lambda'}$, the
  operators $T_{\ell\ell'}$ and $T_{\lambda\lambda'}$ satisfy the following
  exchange algebra
  \begin{equation}
    \label{eq:exchangeTT}
    T_{\ell\ell'}(z_1) \; T_{\lambda\lambda'}(z_2) =
    \mathbf{F}_{\ell\ell'}^{\lambda\lambda'}(z_1/z_2) \;
    T_{\lambda\lambda'}(z_2) \; T_{\ell\ell'}(z_1)
  \end{equation}
  On the surface $\cP_{-\ell,-\ell'}\cap\cP_{-\lambda,-\lambda'}$, the
  operators $S_{\ell\ell'}$ and $S_{\lambda\lambda'}$ satisfy the following
  exchange algebra 
  \begin{equation}
    \label{eq:exchangeSS}
    S_{\ell\ell'}(z_1) \; S_{\lambda\lambda'}(z_2) =
    \mathbf{F}_{\ell\ell'}^{\lambda\lambda'}(z_1/z_2) \;
    S_{\lambda\lambda'}(z_2) \; S_{\ell\ell'}(z_1)
  \end{equation}
  On the surface $\cP_{\ell,\ell'}\cap\cP_{-\lambda,-\lambda'}$, the
  operators $T_{\ell\ell'}$ and $S_{\lambda\lambda'}$ satisfy the following
  exchange algebra 
  \begin{equation}
    \label{eq:exchangeTS}
    T_{\ell\ell'}(z_1) \; S_{\lambda\lambda'}(z_2) = 
    \mathbf{F}_{\ell\ell'}^{\lambda\lambda'}(z_1/z_2) \;
    S_{\lambda\lambda'}(z_2) \; T_{\ell\ell'}(z_1)
  \end{equation}
  The exchange function $\mathbf{F}_{\ell\ell'}^{\lambda\lambda'}(z)$
  is given by
  \begin{equation}
    \label{eq:functstrucTT}
    \mathbf{F}_{\ell\ell'}^{\lambda\lambda'}(z) =
    \frac{F(\ell,z)}{F^*(\ell',z)} \;
    \frac{F^*(\ell',\gamma_{\lambda\lambda'}^{-1}z)}
    {F(\ell,\gamma_{\lambda\lambda'}^{-1}z)} 
  \end{equation}
\end{proposition}
Note that when $c$ is rational, the surface condition implies
the existence of integers $r,s$ such that $p^r=q^s$.\\
Let us remark that eqs. (\ref{eq:exchangeTT}) and (\ref{eq:exchangeSS}) imply the
following compatibility condition:
\begin{equation}
  \label{eq:compat}
  \mathbf{F}_{\ell\ell'}^{\lambda\lambda'}(z) \,
  \mathbf{F}_{\lambda\lambda'}^{\ell\ell'}(z^{-1}) = 1 
\end{equation}
which is trivially satisfied thanks to the properties of the function
$F(\ell,z)$ and the surface conditions. \\
In the particular case $(\lambda,\lambda') = (\ell,\ell')$ which will 
be studied below, formula
(\ref{eq:functstrucTT}) simplifies. Indeed, using the ``inversion'' formula
(\ref{eq:inver}) and the $q^2$-periodicity of the function $F(\ell,z)$,
one gets
\begin{equation}
  \mathbf{F}_{\ell\ell'}^{\ell\ell'}(z) =
  \frac{F(|\ell|,z)}{F(|\ell|,(-p^{\half})^{|\ell|}z)} \;
  \frac{F^*(|\ell'|,q^{2}(-p^{*\half})^{|\ell'|}z)}{F^*(|\ell'|,z)}
\end{equation}

\medskip

Propositions \ref{prop:LT} to \ref{prop:TT} generalize the results
obtained in \cite{ENSLAPP649, LAPTH690} for the particular values
$(\ell,\ell')=(\lambda,\lambda')=(2m+1,1)$. 

\medskip

It is possible to define a semi-classical limit of these exchange algebras.
Indeed, the exchange function $\mathbf{F}_{\ell\ell'}^{\ell\ell'}(z)$
degenerates to 1 when $-p^{\half} = q^k$ ($k \in \ZZ$) and $c$ integer.
Setting $-p^{\half} = q^{k -\beta/\eta}$ for some integer $k \in \ZZ$ and
$\eta = (\ell-\ell')(\ell+\ell'-1)$, a Poisson structure, when $\ell \ne
\ell'$, can now be defined by the following limit\footnotemark:
\begin{equation}
  \Big\{ T_{\ell\ell'}(z_{1}) \, , \, T_{\ell\ell'}(z_{2})\Big\} \equiv
  \lim_{\beta \rightarrow 0} \, \frac{1}{\beta} \, \Big(
  T_{\ell\ell'}(z_{1}) \, T_{\ell\ell'}(z_{2})-T_{\ell\ell'}(z_{2}) \,
  T_{\ell\ell'}(z_{1}) \Big)
  \label{eq:poiss}
\end{equation}
In fact, when computing this limit, one finds that it does not depend on 
$k$, and the result is 
% (with $z_{21} = z_{2}/z_{1}$)
% \begin{eqnarray}
%    \Big\{ T_{\ell\ell'}(z_{1}) \, , \, T_{\ell\ell'}(z_{2})\Big\} &=& 2 \ln
%    q \left[ \frac{z_{21}^2}{1-z_{21}^2} - \frac{z_{21}^{-2}}{1-z_{21}^{-2}} +
%    2\sum_{n=0}^{\infty} \left( - \frac{1}{1-q^{4n}z_{21}^2} +
%    \frac{1}{1-q^{4n+2}z_{21}^2} \right. \right . \nonumber \\
%    && \hspace{10mm} \left. \left. + \; \frac{1}{1-q^{4n}z_{21}^{-2}} -
%    \frac{1}{1-q^{4n+2}z_{21}^{-2}} \right) \right] T_{\ell\ell'}(z_{1}) \,
%    T_{\ell\ell'}(z_{2})
%   \label{eq:resupoiss}
% \end{eqnarray}
\begin{equation}
   \Big\{ T_{\ell\ell'}(z_{1}) \, , \, T_{\ell\ell'}(z_{2})\Big\} =
   h(z_{2}/z_{1}) \, T_{\ell\ell'}(z_{1}) \, T_{\ell\ell'}(z_{2})
  \label{eq:resupoiss}
\end{equation}
where
\begin{equation}
  h(x) = 2 \ln q \left[ \frac{1+x^2}{1-x^2} + 2\sum_{n=0}^{\infty} \left(
  \frac{1}{1-q^{4n+2}x^2} - \frac{1}{1-q^{4n}x^2} + \;
  \frac{1}{1-q^{4n}x^{-2}} - \frac{1}{1-q^{4n+2}x^{-2}} \right) \right]
  \label{eq:foncpoiss}
\end{equation}
This Poisson structure is identical to the classical deformed Virasoro
algebra constructed in \cite{FR}. It is therefore consistent to consider
all these quantum exchange algebras as ``quantum deformed Virasoro
algebras''. 

\footnotetext{When $\ell = \ell'$, the exchange function
$\mathbf{F}_{\ell\ell'}^{\ell\ell'}(z)$ degenerates to 1 for any $p$ and
$q$ and $c = -2/\ell$. The Poisson structure in this case is obtained by 
expanding around this ``critical'' value of $c$, see \cite{FR,ENSLAPP649}.}

\subsection{Example: $c=1$ and $p=q^{3}$}

We consider the case where $c=1$ and the parameters $p$ and $q$ are
related by $p = q^{3}$. The structure constants
$\mathbf{F}_{\ell\ell'}^{\lambda\lambda'}(z)$ depend only on the
congruency classes $\overline{\ell}$ and $\overline{\lambda}$ modulo 4 of
the integers $\ell$ and $\lambda$, the surface conditions implying that the
integers $\ell'$ and $\lambda'$ are given by 
\begin{eqnarray}
\ell'&=&3\ell+4\quad ;\quad \lambda'\ =\ 3\lambda+4\\
\overline\ell' &\equiv& 3\overline\ell \mod 4\quad ;\quad
\overline\lambda'\ \equiv\ 3\overline\lambda \mod 4\,. 
\end{eqnarray}
When $(\overline{\lambda},\overline{\lambda}')=(0,0)$ or
$(\overline{\ell},\overline{\ell}')=(0,0)$, one gets 
$\mathbf{F}_{\ell\ell'}^{\lambda\lambda'}(z)=1$
The non-trivial structure constants
$\mathbf{F}_{\ell\ell'}^{\lambda\lambda'}(z)$ are summarized in the
following tableau ($\Theta(x)\equiv\Theta_{q^4}(x)$):
\begin{displaymath}
  \begin{array}{l|ccc}
  \!\!\begin{array}{lr}
  \qquad 
%   \setminus 
  &\!\!\!\!\!\mbox{\small{$(\overline{\lambda},\overline{\lambda}')$}} \\
\mbox{\small{$(\overline{\ell},\overline{\ell}')$}} 
&\qquad
\end{array}\!\!
    & \begin{array}{c} (1,3)\\ 
   \ \end{array}
   & \begin{array}{c} (2,2)\\ \ \end{array} & \begin{array}{c} (3,1)\\
   \ \end{array} \\
%  	 &&&& \\
   \hline\\
    ( 1,3) &  q^{2} \;
    \dfrac{\Theta(q^{3}z^{2})^{2}\,\Theta(q^{-1}z^{2})^{2}}
    {\Theta(qz^{2})^{4}} & q^{2} \;
    \dfrac{\Theta(q^{2}z^{2})^{4}\,\Theta(q^{-1}z^{2})^{4}}
    {\Theta(z^{2})^{4}\,\Theta(qz^{2})^{4}} & q^{2} \;
    \dfrac{\Theta(q^{2}z^{2})^{2}\,\Theta(q^{-2}z^{2})^{2}}
    {\Theta(z^{2})^{4}} \\
    &&& \\
    (2,2) &  q^{2} \;
    \dfrac{\Theta(z^{2})^{4}\,\Theta(q^{3}z^{2})^{4}}
    {\Theta(q^{2}z^{2})^{4}\,\Theta(qz^{2})^{4}} & q^{4}
    \;
    \dfrac{\Theta(q^{3}z^{2})^{4}\,\Theta(q^{-1}z^{2})^{4}}
    {\Theta(qz^{2})^{8}} & q^{2} \;
    \dfrac{\Theta(q^{2}z^{2})^{4}\,\Theta(q^{-1}z^{2})^{4}}
    {\Theta(z^{2})^{4}\,\Theta(qz^{2})^{4}} \\
    &&& \\
    (3,1) & q^{-2} \; \dfrac{\Theta(z^{2})^{4}}
    {\Theta(q^{2}z^{2})^{2}\,\Theta(q^{-2}z^{2})^{2}} &
    q^{2} \;
    \dfrac{\Theta(z^{2})^{4}\,\Theta(q^{3}z^{2})^{4}}
    {\Theta(q^{2}z^{2})^{4}\,\Theta(qz^{2})^{4}} & q^{2}
    \;
    \dfrac{\Theta(q^{3}z^{2})^{2}\,\Theta(q^{-1}z^{2})^{2}}
    {\Theta(qz^{2})^{4}} \\
  \end{array}
\end{displaymath}

\subsection{Riemann--Hilbert splitting}

In order to define the algebra generated by $T_{\ell\ell'}(z)$ and its
exchange relations, one has to introduce the modes of the generators
$T_{\ell\ell'}(z)$, i.e. $T_{\ell\ell'}(z) = \sum_{n \in \ZZ}
T_{\ell\ell'}(n) \, z^{-n}$, and use the exchange relations
(\ref{eq:exchangeTT}) (with $\lambda,\lambda' = \ell,\ell'$) as ordering
relations among the modes $T_{\ell\ell'}(n)$ (see e.g. \cite{FIJKMY95}).
For this purpose, we need to prepare a Riemann--Hilbert splitting of the
exchange function in (\ref{eq:functstrucTT}), i.e. to factorize the
exchange function into a function analytic around $z = 0$ and a function
analytic around $z^{-1}=0$:
\begin{eqnarray}
  \mathbf{F}_{\ell\ell'}^{\ell\ell'}(z) &=& \boldsymbol{\varphi}_+(z)
  \; \boldsymbol{\varphi}_-^{-1}(z^{-1})\\
  \boldsymbol{\varphi}_{\pm}(z) &=& 1+O(z^{\pm1})
\end{eqnarray}
in the neighborhood of a circle $C$ of radius $R$. $\boldsymbol{\varphi}_+$ and
$\boldsymbol{\varphi}_-$ are respectively analytic for $|z|<R$ and $|z|>R$.
It may be possible to choose $\boldsymbol{\varphi}_+ = \boldsymbol{\varphi}_- =
\varphi$ in the sense of analytic continuation. Indeed, in our case, the exchange
relation (\ref{eq:exchangeTT}) reads, after the Riemann--Hilbert splitting,
as
\begin{equation}
  \varphi_{\ell\ell'}(z_2/z_1) \; T_{\ell\ell'}(z_1) \;
  T_{\ell\ell'}(z_2) = \varphi_{\ell\ell'}(z_1/z_2) \;
  T_{\ell\ell'}(z_2) \; T_{\ell\ell'}(z_1)
\end{equation}
where 
% ANCIENNE FORMULE
% \begin{eqnarray}
%   && \varphi_{\ell\ell'}(x) = \Bigg[ \frac{1}{(1-x^2)^{\ell+\ell'-2}} \;
%   \prod_{s=1}^{\ell-1} \frac{(p^s x^2;q^4)_\infty \, (q^2 p^{-s}
%   x^2;q^4)_\infty} {(q^4 p^{-s} x^2;q^4)_\infty \, (q^2 p^s
%   x^2;q^4)_\infty} \; \nonumber \\
%   && \hspace*{116pt} \times \; \prod_{s=1}^{\ell'-1} \frac{({p^*}^{-s}
%   x^2;q^4)_\infty \, (q^2 {p^*}^s x^2;q^4)_\infty} {(q^4 {p^*}^s
%   x^2;q^4)_\infty \, (q^2 {p^*}^{-s} x^2;q^4)_\infty} \Bigg]^2
% \end{eqnarray}
% NOUVELLE FORMULE   \; \nonumber \\ && \hspace*{116pt}
\begin{eqnarray*}
 \varphi_{\ell\ell'}(x) &=& \Bigg[ \frac{1}{(1-x^2)^{|\ell|-|\ell'|}} \;
  \prod_{s=1}^{|\ell|-1} \frac{(p^s x^2;q^4)_\infty \, (q^2 p^{-s}
  x^2;q^4)_\infty} {(q^4 p^{-s} x^2;q^4)_\infty \, (q^2 p^s
  x^2;q^4)_\infty} 
  \\
  &&\hspace{16ex} \prod_{s=1}^{|\ell'|-1} \frac{(q^4
  {p^*}^{-s} x^2;q^4)_\infty \, (q^2 {p^*}^s x^2;q^4)_\infty} {({p^*}^s
  x^2;q^4)_\infty \, (q^2 {p^*}^{-s} x^2;q^4)_\infty}
  \Bigg]^2\qquad
\end{eqnarray*}
This choice of analyticity properties for the exchange functions guarantees
the existence of a consistent normal ordering procedure based on the
reordering by increasing values of the mode indices. This allows one in
turn to define a suitable Poincar\'e--Birkhoff--Witt basis for the algebra.

Once this algebra is well-defined, one can wonder whether it admits a 
central extension term. From the 
results of \cite{AKOS96}, one already knows that the central extension will 
be different from the one computed in this reference. Indeed, it was proved
there that this central extension uniquely determine the exchange
function, which is not of the type $\phi_{\ell\ell'}(x)$. However, in 
our case, other kinds of term may arise. Hence, this question remains open. 

\section{Vertex operators\label{sect:VO}}

An interesting interpretation of the surfaces arises when considering
vertex operator representations \cite{FIJKMY95,Junichi04} of the
elliptic algebra ${\cA}_{q,p}(\widehat{sl}_{2})$ at $c=1$. These surfaces
are related to coincident singularities in the Riemann--Hilbert splitted
exchange relations of the vertex operators.

\subsection{Level one vertex operators}

The well-known Verma modules and irreducible highest weight modules of the
affine Lie algebras are expected to have deformations to the elliptic case.
Denoting by $V(\Lambda_{i})$ and $V(\Lambda_{1-i})$ ($i=0,1$) the
${\cA}_{q,p}(\widehat{sl}_{2})$-modules corresponding to the level one
irreducible highest weight modules of ${\cU}_{q}(\widehat{sl}_{2})$, one
introduces the level one vertex operators defined as intertwiners between
${\cA}_{q,p}(\widehat{sl}_{2})$-modules (see \cite{FIJKMY95,Junichi04}):
\begin{eqnarray}
  && \Phi^{(i)}(z) : V(\Lambda_{i}) \longrightarrow V(\Lambda_{1-i})
  \otimes V_{z} \\
  && \Psi^{*(i)}(z) : V_{z} \otimes V(\Lambda_{i}) \longrightarrow
  V(\Lambda_{1-i})
\end{eqnarray}
where $V_{z} = \CC[z,z^{-1}] \otimes (\CC v_{+} \oplus \CC v_{-})$ is the
spin 1/2 evaluation module.

One assumes existence and uniqueness of these operators in the elliptic
case. $\Phi$ and $\Psi^{*}$ are called type I and type II vertex operators.
\\
We will use the following decomposition
\begin{eqnarray}
  \Phi^{(i)}(z) &=& \Phi^{(i)}_{+}(z) \otimes v_{+} + \Phi^{(i)}_{-}(z) 
  \otimes v_{-} \quad
  \mbox{with} \quad 
  \Phi^{(i)}_{\eps}(-z) = -(-1)^{i}\,\eps\,\Phi^{(i)}_{\eps}(z) 
%   \Phi^{(i)}_{\eps}(z) = \sum_{n \equiv (1 \pm
%   \eps)/2} \Phi^{(i)}_{\eps,n} z^{-n}
  \label{eq:modesphi}
\\
\Psi^{*(i)}(z) &=& \Psi^{*(i)}_{+}(z) + \Psi^{*(i)}_{-}(z)
  \quad \mbox{with} \quad 
  \Psi^{*(i)}_{\eps}(-z) = (-1)^{i}\,\eps\,\Psi^{*(i)}_{\eps}(z) 
%   \Psi^{*(i)}_{\eps}(z) = \sum_{n \equiv (1 \mp
%   \eps)/2} \Psi^{*(i)}_{\eps,n} z^{-n}
  \label{eq:modespsi}
\end{eqnarray}
It is conjectured\footnotemark~that the commutation relations between the
vertex operators are given by 
\footnotetext{It was shown in \cite{JKOS} that such relations are valid
when $p$ is infinitesimally small. Their validity for finite $p$ remains an
open question.}
\begin{eqnarray}
  \label{eq:commrel1}
  && \Phi^{(1-i)}_{\eps_{2}}(z_2) \, \Phi^{(i)}_{\eps_{1}}(z_1) =
  \sum_{\eps'_{1},\eps'_{2}=\pm}
  \cW_{\eps_{1}\eps_{2}}^{\eps'_{1}\eps'_{2}} (z_1/z_2)
  \, \Phi^{(1-i)}_{\eps'_{1}}(z_1) \, \Phi^{(i)}_{\eps'_{2}}(z_2) \\
  \label{eq:commrel2}
   && \Psi^{*(1-i)}_{\eps_{1}}(z_1) \, \Psi^{*(i)}_{\eps_{2}}(z_2) = -
   \sum_{\eps'_{1},\eps'_{2}=\pm} \Psi^{*(1-i)}_{\eps'_{2}}(z_2) \,
   \Psi^{*(i)}_{\eps'_{1}}(z_1) \,
   {\cW^*}_{\eps_{1}\eps_{2}}^{\eps'_{1}\eps'_{2}} (z_1/z_2) \\
   \label{eq:commrel3}
  && \Phi^{(1-i)}_{\eps_{1}}(z_1) \, \Psi^{*(i)}_{\eps_{2}}(z_2) = \tau
  (z_1/z_2) \, \Psi^{*(1-i)}_{\eps_{2}}(z_2) \, 
  \Phi^{(i)}_{\eps_{1}}(z_1)
\end{eqnarray}

In terms of the vertex operators, the generators $L(z)$ are given by
\begin{equation}
  L_{\eps\eps'}(z) = \kappa \Psi^{*}_{\eps'}(q^{-\half}z)
  \Phi_{\eps}(z)
  \label{eq:Lvertex}
\end{equation}
while their inverse read
\begin{equation}
  L^{-1}_{\eps\eps'}(z) = \kappa^{-1}
  \Psi^{*}_{-\eps}(-q^{-\three}z) \Phi_{-\eps'}(-q^{-1}z)\,.
  \label{eq:Lvertexinv}
\end{equation}

% Let us remark that one can gather the operators $\Phi^{(i)}_{\eps}(z)$ and 
% $\Psi^{*(i)}_{\eps}(z)$ into column and row vectors that we will
% loosely write $\Phi$ and $\Psi^{\dag}$:
% \begin{equation}
% \Phi(z)=\left(\begin{array}{c} \Phi_{+}(z) \\ \Phi_{-}(z)
% \end{array}\right)
% \quad\mbox{and}\quad
% \Psi^\dag(z) = \Big( \Psi_{+}^{*}(z)\,,\ \Psi_{-}^{*}(z) \Big)\,.
% \end{equation}
% Then, using auxiliary spaces, the above relations can be recasted as
% \begin{eqnarray}
%   \label{eq:commrelaux1}
%   && \Phi_{2}(z_2) \, \Phi_{1}(z_1) = R_{12}(z_1/z_2) \,
%   \Phi_{1}(z_1) \, \Phi_{2}(z_2) \\
%   \label{eq:commrelaux2}
%   && \Psi^{\dag}_{1}(z_1) \, \Psi^{\dag}_{2}(z_2) = -
%   \Psi^{\dag}_{2}(z_2) \, \Psi^{\dag}_{1}(z_1) \,
%   {R^*}_{12}(z_1/z_2) \\
%   \label{eq:commrelaux3}
%   && \Phi_{1}(z_1) \, \Psi^{\dag}_{2}(z_2) = \tau(z_1/z_2)
%   \, \Psi^{\dag}_{2}(z_2) \, \Phi_{1}(z_1)
% \end{eqnarray}
% \begin{eqnarray}
%   L(z) &=& \kappa\ i_{12}\ L_{12}(z) 
%   = \kappa\ i_{12}\ \Psi^{*}_{2}(q^{-\half}z) \Phi_{1}(z)\\
%   L^{-1}(z) &=& \kappa^{-1}\ i_{12}\ L_{21}(-q^{-1}z) 
%   = \kappa^{-1}\ i_{12}\ \Psi^{*}_{1}(-q^{-\three}z) 
%   \Phi_{2}(-q^{-1}z)
% \end{eqnarray}
% where $i_{12}$ is the identification of the tensor product of a row
% and a column vector with a matrix.

\subsection{Connection with the surface conditions\label{sect:residu}}

{From} the relations
\begin{align}
  & \rho(z) \big( a(z) \pm d(z) \big) = z^{-1}
  \frac{\alpha^{\pm}(z^{-1})}{\alpha^{\pm}(z)} \;, && \qquad
  \alpha^{\pm}(z) = \frac{(\pm p^{\half}qz;p)_{\infty}}{(\pm
  p^{\half}q^{-1}z;p)_{\infty}} \; \xi(z^2;p,q^4)^{-1} \\
  & \rho(z) \big( b(z) \pm c(z) \big) =
  \frac{\beta^{\pm}(z^{-1})}{\beta^{\pm}(z)} \;, && \qquad
  \beta^{\pm}(z) = \frac{(\mp qz;p)_{\infty}}{(\mp
  pq^{-1}z;p)_{\infty}} \; \xi(z^2;p,q^4)^{-1}
\end{align}
we have the following Riemann--Hilbert splittings \cite{FIJKMY95,Junichi04}:
\begin{eqnarray}
 &&  z_{1}^{-1} \alpha^{\pm}(z_{2}/z_{1}) \; \Big(
   \Phi_{+}(z_{1}) \Phi_{+}(z_{2}) \pm \Phi_{-}(z_{1})
   \Phi_{-}(z_{2}) \Big) \nonumber\\
&&\qquad\qquad\qquad =z_{2}^{-1}
   \alpha^{\pm}(z_{1}/z_{2}) \; \Big( \Phi_{+}(z_{2})
   \Phi_{+}(z_{1}) \pm \Phi_{-}(z_{2}) \Phi_{-}(z_{1}) \Big)\qquad\qquad
  \label{eq:RHsplitphiA}
\\
  && \beta^{+}(z_{2}/z_{1}) \; \Big( \Phi_{+}(z_{1})
   \Phi_{-}(z_{2}) + \Phi_{-}(z_{1}) \Phi_{+}(z_{2}) \Big) \nonumber\\
 &&\qquad\qquad\qquad =\beta^{+}(z_{1}/z_{2}) \; \Big( \Phi_{+}(z_{2})
   \Phi_{-}(z_{1}) + \Phi_{-}(z_{2}) \Phi_{+}(z_{1}) \Big)\qquad\qquad
\label{eq:RHsplitphiB} 
\end{eqnarray}
and
\begin{eqnarray}
   && \frac{z_{1}^{-1}}{1-z_{2}^2/z_{1}^2} \;
   \frac{1}{\alpha^{*\pm}(z_{2}/z_{1})} \; \Big(
   \Psi^{*}_{+}(z_{1}) \Psi^{*}_{+}(z_{2}) \pm
   \Psi^{*}_{-}(z_{1}) \Psi^{*}_{-}(z_{2}) \Big) \nonumber \\
   && \qquad\qquad \;\;=\;\;
   \frac{z_{2}^{-1}}{1-z_{1}^2/z_{2}^2} \; 
   \frac{1}{\alpha^{*\pm}(z_{1}/z_{2})} \; \Big(
   \Psi^{*}_{+}(z_{2}) \Psi^{*}_{+}(z_{1}) \pm
   \Psi^{*}_{-}(z_{2}) \Psi^{*}_{-}(z_{1}) \Big) \qquad 
   \label{eq:RHsplitpsiA} \\
   && \nonumber\\
   && \frac{(1+q^{-1}z_{2}/z_{1})(1+qz_{2}/z_{1})}
   {1-z_{2}^2/z_{1}^2} \;
   \frac{1}{\beta^{*+}(z_{2}/z_{1})} \; \Big(
   \Psi^{*}_{+}(z_{1}) \Psi^{*}_{-}(z_{2}) +
   \Psi^{*}_{-}(z_{1}) \Psi^{*}_{+}(z_{2}) \Big) \nonumber \\
   &&  \;\;=\;\;
   \frac{(1+q^{-1}z_{1}/z_{2})(1+qz_{1}/z_{2})}
   {1-z_{1}^2/z_{2}^2} \;
   \frac{1}{\beta^{*+}(z_{1}/z_{2})} \; \Big(
   \Psi^{*}_{+}(z_{2}) \Psi^{*}_{-}(z_{1}) +
   \Psi^{*}_{-}(z_{2}) \Psi^{*}_{+}(z_{1}) \Big) \qquad
   \label{eq:RHsplitpsiB} 
\end{eqnarray}
\begin{equation}
  \label{eq:RHsplitphipsi}
    z_{1} \; \frac{(qz_{2}^{2}/z_{1}^{2};q^4)_{\infty}}
    {(q^3z_{2}^{2}/z_{1}^{2};q^4)_{\infty}} \;
    \Phi_{\eps_{1}}(z_{1}) \Psi^{*}_{\eps_{2}}(z_{2}) =
    z_{2} \; \frac{(qz_{1}^{2}/z_{2}^{2};q^4)_{\infty}}
    {(q^3z_{1}^{2}/z_{2}^{2};q^4)_{\infty}} \;
    \Psi^{*}_{\eps_{2}}(z_{2}) \Phi_{\eps_{1}}(z_{1})
\end{equation}
As before, $\alpha^{*\pm}(z)$ and $\beta^{*\pm}(z)$ correspond to
$\alpha^{\pm}(z)$ and $\beta^{\pm}(z)$ with $p \to p^* = pq^{-2c}$.

\medskip

Equations (\ref{eq:RHsplitphiA})--(\ref{eq:RHsplitpsiB}) all have the form
\begin{equation}
  c(z_{2}/z_{1}) \, {\cal O}(z_{1}) \, {\cal O}(z_{2}) = 
  c(z_{1}/z_{2}) \, {\cal O}(z_{2}) \, {\cal O}(z_{1})
  \label{eq:coocoo}
\end{equation}
Our interpretation of the surface conditions will be based on the following
reading of any exchange relation (\ref{eq:coocoo}). Suppose that the
coefficient $c(z_{1}/z_{2})$ of the r.h.s. of (\ref{eq:coocoo}) exhibits a
zero at some position $z_{1}/z_{2} = z_{0}$ and that the coefficient
$c(z_{2}/z_{1})$ of the l.h.s. of (\ref{eq:coocoo}) does not have any zero
nor pole at $z_{1}/z_{2} = z_{0}$. Thus, for the equality to hold, the
operator ${\cal O}(z_{2}) \, {\cal O}(z_{1})$ on the r.h.s. must have a
pole at $z_{1}/z_{2} = z_{0}$ and the operator ${\cal O}(z_{1}) \, {\cal
O}(z_{2})\big\vert_{z_{1}/z_{2} = z_{0}}$ on the l.h.s. is then interpreted
as a \textit{residue} operator. In the same way, when the coefficient
$c(z_{2}/z_{1})$ of the l.h.s. of (\ref{eq:coocoo}) has a pole at some
position $z_{1}/z_{2} = z'_{0}$ for which the coefficient $c(z_{1}/z_{2})$
of the r.h.s. is regular, the bilinear ${\cal O}(z_{2}) \, {\cal O}(z_{1})$
on the r.h.s. has a pole at the same location, and ${\cal O}(z_{1}) \,
{\cal O}(z_{2})\big\vert_{z_{1}/z_{2} = z'_{0}}$ on the l.h.s. is
interpreted as a residue operator. In both cases, it may also occur that
${\cal O}(z_{2}) \, {\cal O}(z_{1})$ be regular and ${\cal O}(z_{1}) \,
{\cal O}(z_{2})\big\vert_{z_{1}/z_{2} = z_{0} \textrm{ or } z'_{0}}$
degenerate to zero, still interpreted as a ``residue'' of a regular
operator.

\medskip

The aim of this paragraph is to propose an interpretation of the operators
$T_{\ell\ell'}$ and $S_{\ell\ell'}$ as bilinear of such residue operators
in the context of the level one vertex operator representation. For such a
purpose, let us define the following operator
\begin{equation}
  \label{eq:defT12}
  \cT_{\bar\ell\bar\ell'}(z_{1},z_{2}) \equiv \tr \Big( L(z_{1})^{-1} \;
  \sigma_{x}^{\ell} \; L(z_{2}) \; \sigma_{x}^{\ell'} \Big)
  \qquad \ell\,,\ \ell'\in\ZZ
\end{equation}
Because of the property of Pauli matrices, it depends on $\ell$ and
$\ell'$ only modulo 2, hence the notation $\bar\ell$ and $\bar\ell'$.
This operator is formally related to $T_{\ell\ell'}$:
\begin{equation}
  T_{\ell\ell'}(z)=\cT_{\bar\ell\bar\ell'}(\gamma_{\ell\ell'}z,z) \,.
\end{equation}
It is this formal relation which we shall now clarify. \\
In terms of type I and type II vertex operators, one gets for $\ell$ and
$\ell'$ even,
\begin{eqnarray}
  \cT_{00}(z_{1},z_{2}) &=&
  \Psi^{*}_{+}(-q^{-\three}z_{1}) \; \Phi_{+}(-q^{-1}z_{1}) \;
  \Psi^{*}_{-}(q^{-\half}z_{2}) \; \Phi_{-}(z_{2}) \nonumber \\
  &&+\  \Psi^{*}_{+}(-q^{-\three}z_{1}) \; \Phi_{-}(-q^{-1}z_{1})
  \; \Psi^{*}_{-}(q^{-\half}z_{2}) \; \Phi_{+}(z_{2}) \nonumber
  \\
  &&+\  \Psi^{*}_{-}(-q^{-\three}z_{1}) \; \Phi_{+}(-q^{-1}z_{1})
  \; \Psi^{*}_{+}(q^{-\half}z_{2}) \; \Phi_{-}(z_{2}) \nonumber
  \\
  &&+\  \Psi^{*}_{-}(-q^{-\three}z_{1}) \; \Phi_{-}(-q^{-1}z_{1})
  \; \Psi^{*}_{+}(q^{-\half}z_{2}) \; \Phi_{+}(z_{2})
  \label{eq:T00vertex}
\end{eqnarray}
Then using eq. (\ref{eq:RHsplitphipsi}), one obtains
\begin{eqnarray}
 && \frac{1}{q^{-\half}z_{2}} \;
  \frac{(q^2z_{2}^{2}/z_{1}^{2};q^4)_{\infty}}
  {(q^4z_{2}^{2}/z_{1}^{2};q^4)_{\infty}} \;
  \cT_{00}(z_{1},z_{2}) = \frac{1}{q^{-1}z_{1}} \;
  \frac{(z_{1}^{2}/z_{2}^{2};q^4)_{\infty}}
  {(q^2z_{1}^{2}/z_{2}^{2};q^4)_{\infty}} \; \nonumber \\
 & &\qquad\qquad\qquad\qquad\times \Big( \Psi^{*}_{+}(-q^{-\three}z_{1}) \;
  \Psi^{*}_{-}(q^{-\half}z_{2}) + \Psi^{*}_{-}(-q^{-\three}z_{1})
  \; \Psi^{*}_{+}(q^{-\half}z_{2}) \Big) \nonumber \\
 & &\qquad\qquad\qquad\qquad\times \Big( \Phi_{+}(-q^{-1}z_{1}) \; \Phi_{-}(z_{2}) +
  \Phi_{-}(-q^{-1}z_{1}) \; \Phi_{+}(z_{2}) \Big)
  \label{eq:T00vertexret}
\end{eqnarray}
We observe that $\cT_{00}(z_{1},z_{2})$ is expressed in terms of the
bilinear of vertex operators $\Phi$ and $\Psi$, which appear in the l.h.s.
of the Riemann--Hilbert splittings
(\ref{eq:RHsplitphiA})--(\ref{eq:RHsplitpsiB}). Our previous discussion on
the meaning of eq. (\ref{eq:coocoo}) indicates that we must now carefully
analyze the pole and zero structure of the coefficients of
(\ref{eq:RHsplitphiA})--(\ref{eq:RHsplitpsiB}). The operators
$T_{\ell\ell'}$ and the associated surface conditions will in fact appear
from the requirement of simultaneous singularities in both bilinears
$\Phi\Phi$ and $\Psi\Psi$. \\
We must first of all refocus the singularity analysis on modified exchange 
relations (\ref{eq:RHsplitphiA})--(\ref{eq:RHsplitpsiB}).
The structure functions $\alpha^\pm(z)$ and $\beta^\pm(z)$ all contain the
same normalization function $\xi(z^2;p,q^4)$. As a working hypothesis,
we will not consider the poles arising from this function. A
heuristic argument for such a point of view can be sketched as follows.
The function $\xi(z^2;p,q^4)$ arises from the
normalization coefficient of the elliptic $R$ matrix and can be viewed as
the contribution of a $U(1)$ current. However, this $U(1)$ current must not
be confused with the (elliptic analogue of the) $gl(1)$ current completing
$sl(2)$ into $gl(2)$. Indeed, in our context, this latter current has
vanishing commutation relations with the ${\cA}_{q,p}(\widehat{gl}_{2})$
currents (including itself). It is possible to use the $U(1)$ current in a
redefinition of the vertex operators to eliminate the $\xi$ factor. This
redefinition implies for the vertex operators $\Phi$ and $\Psi$ of
${\cA}_{q,p}(\widehat{sl}_{2})$ an extension of the Verma module to
incorporate this extra $U(1)$ current. For the sake of simplicity, we shall
not implement this extension here, but its existence shows that one can
consistently restrict the analysis of the poles and zeroes in the exchange
relations (\ref{eq:RHsplitphiA})--(\ref{eq:RHsplitpsiB}) to the part of the
$\alpha^\pm(z)$ and $\beta^\pm(z)$ functions factoring the $\xi$ function.
From now on, this restriction will be implicit. \\
Consider first equation (\ref{eq:RHsplitphiB}):
\begin{eqnarray}
  \label{eq:splitphiB}
  && \beta^{+}(-qz_{2}/z_{1}) \; \Big( \Phi_{+}(-q^{-1}z_{1})
  \; \Phi_{-}(z_{2}) + \Phi_{-}(-q^{-1}z_{1}) \;
  \Phi_{+}(z_{2}) \Big) \;\;=\;\; \nonumber \\
  && \qquad \qquad \qquad \beta^{+}(-q^{-1}z_{1}/z_{2}) \; \Big(
  \Phi_{+}(z_{2}) \; \Phi_{-}(-q^{-1}z_{1}) + \Phi_{-}(z_{2})
  \; \Phi_{+}(-q^{-1}z_{1}) \Big)\qquad\qquad
\end{eqnarray}
At the values $z_{1}/z_{2} = (-p^{\half})^{-\ell}$ where $\ell \in
2\ZZ_{\ge 0}$, the function $\beta^{+}$ of the r.h.s. of eq.
(\ref{eq:splitphiB}) has a zero while the function $\beta^{+}$ of the
l.h.s. has neither poles nor zeroes. Similarly, at the values $z_{1}/z_{2}
= (-p^{\half})^{-\ell}$ where $\ell \in 2\ZZ_{<0}$, the function
$\beta^{+}$ of the l.h.s. of eq. (\ref{eq:splitphiB}) has a pole, while the
function $\beta^{+}$ of the l.h.s. is regular. Hence the operator
$\Phi_{+}(z_{2}) \; \Phi_{-}(-q^{-1}z_{1}) + \Phi_{-}(z_{2}) \;
\Phi_{+}(-q^{-1}z_{1})$ must have a pole located at $z_{1}/z_{2} =
(-p^{\half})^{-\ell}$, $\ell \in 2\ZZ$, while $\Phi_{+}(-q^{-1}z_{1}) \;
\Phi_{-}(z_{2}) + \Phi_{-}(-q^{-1}z_{1}) \; \Phi_{+}(z_{2})$ is interpreted
as a residue operator. \\
In the same way, from equation (\ref{eq:RHsplitpsiB}), one has
\begin{eqnarray}
  \label{eq:splitpsiB}
   && \frac{(1-z_{2}/z_{1})(1-q^2z_{2}/z_{1})}
   {(1-q^2z_{2}^2/z_{1}^2)\;\beta^{*+}(-qz_{2}/z_{1})} \;
    \Big(
   \Psi^{*}_{+}(-q^{-\three}z_{1}) \Psi^{*}_{-}(q^{-\half}z_{2})
   + \Psi^{*}_{-}(-q^{-\three}z_{1})
   \Psi^{*}_{+}(q^{-\half}z_{2}) \Big) \;\;=\;\; \nonumber \\[2.1ex]
   && \frac{(1-z_{1}/z_{2})(1-q^{-2}z_{1}/z_{2})}
   {(1-q^{-2}z_{1}^2/z_{2}^2)\;\beta^{*+}(-q^{-1}z_{1}/z_{2})} \;
   \Big(
   \Psi^{*}_{+}(q^{-\half}z_{2}) \Psi^{*}_{-}(-q^{-\three}z_{1})
   + \Psi^{*}_{-}(q^{-\half}z_{2})
   \Psi^{*}_{+}(-q^{-\three}z_{1}) \Big) \nonumber \\
   &&
\end{eqnarray}
At the values $z_{1}/z_{2} = q^2 (-p^{*\half})^{-\ell'}$ where $\ell' \in
2\ZZ_{\ge 0}$, the function $1/\beta^{*+}$ of the r.h.s. of eq.
(\ref{eq:splitpsiB}) has a zero while the function $1/\beta^{*+}$ of the
l.h.s. has no poles nor zeroes. Similarly, at the values $z_{1}/z_{2} =
q^2 (-p^{*\half})^{-\ell'}$ where $\ell' \in 2\ZZ_{<0}$, the function
$1/\beta^{*+}$ of the r.h.s. of eq. (\ref{eq:splitpsiB}) has a pole while
the function $1/\beta^{*+}$ of the l.h.s. is regular. It follows that in
both cases the quantity $\Psi^{*}_{+}(q^{-\half}z_{2})
\Psi^{*}_{-}(-q^{-\three}z_{1}) + \Psi^{*}_{-}(q^{-\half}z_{2})
\Psi^{*}_{+}(-q^{-\three}z_{1})$ must have a pole at $z_{1}/z_{2} = q^2
(-p^{*\half})^{-\ell'}$, $\ell' \in 2\ZZ$, while
$\Psi^{*}_{+}(-q^{-\three}z_{1}) \Psi^{*}_{-}(q^{-\half}z_{2}) +
\Psi^{*}_{-}(-q^{-\three}z_{1}) \Psi^{*}_{+}(q^{-\half}z_{2})$ is
interpreted as a residue operator. \\
Gathering these analyses, one concludes that the operator
$T_{\ell\ell'}(z)$ has to be interpreted as a product of residue operators
from expression (\ref{eq:T00vertex}) when both conditions $z_{1}/z_{2} =
(-p^{\half})^{-\ell}$ and $z_{1}/z_{2} = q^2 (-p^{*\half})^{-\ell'}$, with
$\ell$ and $\ell'$ even, are fulfilled.

\medskip

Consider now the case where $\ell$ and $\ell'$ are odd. In terms of type I
and type II vertex operators, one gets
\begin{eqnarray}
  \cT_{11}(z_{1},z_{2}) &=&
  \Psi^{*}_{-}(-q^{-\three}z_{1}) \; \Phi_{-}(-q^{-1}z_{1}) \;
  \Psi^{*}_{-}(q^{-\half}z_{2}) \; \Phi_{-}(z_{2}) \nonumber \\
  &&+\  \Psi^{*}_{+}(-q^{-\three}z_{1}) \; \Phi_{+}(-q^{-1}z_{1})
  \; \Psi^{*}_{+}(q^{-\half}z_{2}) \; \Phi_{+}(z_{2}) \nonumber
  \\
  &&+\  \Psi^{*}_{-}(-q^{-\three}z_{1}) \; \Phi_{+}(-q^{-1}z_{1})
  \; \Psi^{*}_{-}(q^{-\half}z_{2}) \; \Phi_{+}(z_{2}) \nonumber
  \\
  &&+\  \Psi^{*}_{+}(-q^{-\three}z_{1}) \; \Phi_{-}(-q^{-1}z_{1})
  \; \Psi^{*}_{+}(q^{-\half}z_{2}) \; \Phi_{-}(z_{2})
  \label{eq:T11vertex}
\end{eqnarray}
and using again (\ref{eq:RHsplitphipsi}), one obtains
\begin{eqnarray}
&&  \frac{1}{q^{-\half}z_{2}} \;
  \frac{(q^2z_{2}^{2}/z_{1}^{2};q^4)_{\infty}}
  {(q^4z_{2}^{2}/z_{1}^{2};q^4)_{\infty}} \;
  \cT_{11}(z_{1},z_{2}) = \frac{1}{q^{-1}z_{1}} \;
  \frac{(z_{1}^{2}/z_{2}^{2};q^4)_{\infty}}
  {(q^2z_{1}^{2}/z_{2}^{2};q^4)_{\infty}} \; \nonumber \\
  &&\qquad\qquad\qquad\qquad\times \Big( \Psi^{*}_{+}(-q^{-\three}z_{1}) \;
  \Psi^{*}_{+}(q^{-\half}z_{2}) + \Psi^{*}_{-}(-q^{-\three}z_{1})
  \; \Psi^{*}_{-}(q^{-\half}z_{2}) \Big) \nonumber \\
  &&\qquad\qquad\qquad\qquad\times \Big( \Phi_{+}(-q^{-1}z_{1}) \; \Phi_{+}(z_{2}) +
  \Phi_{-}(-q^{-1}z_{1}) \; \Phi_{-}(z_{2}) \Big)\qquad\qquad
  \label{eq:T11vertexret}
\end{eqnarray}
{From} equation (\ref{eq:RHsplitphiA}), one has
\begin{eqnarray}
  && \frac{q}{z_{1}} \; \alpha^{+}(-qz_{2}/z_{1}) \; \Big(
  \Phi_{+}(-q^{-1}z_{1}) \Phi_{+}(z_{2}) +
  \Phi_{-}(-q^{-1}z_{1}) \Phi_{-}(z_{2}) \Big) \;\;=\;\; \nonumber
  \\
  && \qquad \qquad \qquad \frac{1}{z_{2}} \;
  \alpha^{+}(-q^{-1}z_{1}/z_{2}) \; \Big( \Phi_{+}(z_{2})
  \Phi_{+}(-q^{-1}z_{1}) + \Phi_{-}(z_{2})
  \Phi_{-}(-q^{-1}z_{1}) \Big)\qquad\qquad
  \label{eq:splitphiA}
\end{eqnarray}
At the values $z_{1}/z_{2} = (-p^{\half})^{-\ell}$ where $\ell \in
\ZZ_{\ge 0}$, $\ell$ odd, the function $\alpha^{+}$ of the r.h.s. of eq.
(\ref{eq:splitphiA}) has a zero while the function $\alpha^{+}$ of the
l.h.s. has neither poles nor zeroes. In addition, at the values
$z_{1}/z_{2} = (-p^{\half})^{-\ell}$ where $\ell \in \ZZ_{<0}$, $\ell$ odd,
the function $\alpha^{+}$ of the l.h.s. of eq. (\ref{eq:splitphiA}) has a
pole while the function $\alpha^{+}$ of the r.h.s. is regular. Hence the
operator $\Phi_{+}(z_{2}) \; \Phi_{+}(-q^{-1}z_{1}) + \Phi_{-}(z_{2}) \;
\Phi_{-}(-q^{-1}z_{1})$ has a pole at $z_{1}/z_{2} = 
(-p^{\half})^{-\ell}$, $\ell+1 \in 2\ZZ$,
while $\Phi_{+}(-q^{-1}z_{1}) \; \Phi_{+}(z_{2}) + \Phi_{-}(-q^{-1}z_{1})
\; \Phi_{-}(z_{2})$ is interpreted as a residue operator. \\
Similarly, from equation (\ref{eq:RHsplitpsiA}), one has
\begin{eqnarray}
  \label{eq:splitpsiA}
   && \frac{z_{1}^{-1}} {1-q^2z_{2}^2/z_{1}^2} \;
   \frac{1}{\alpha^{*+}(-qz_{2}/z_{1})} \; \Big(
   \Psi^{*}_{+}(-q^{-\three}z_{1}) \Psi^{*}_{+}(q^{-\half}z_{2}) +
   \Psi^{*}_{-}(-q^{-\three}z_{1}) \Psi^{*}_{-}(q^{-\half}z_{2})
   \Big) \;\;=\;\; \nonumber \\[2.1ex]
   && \frac{-q^{-1}z_{2}^{-1}} {1-q^{-2}z_{1}^2/z_{2}^2} \;
   \frac{1}{\alpha^{*+}(-q^{-1}z_{1}/z_{2})} \; \Big(
   \Psi^{*}_{+}(q^{-\half}z_{2}) \Psi^{*}_{+}(-q^{-\three}z_{1}) +
   \Psi^{*}_{-}(q^{-\half}z_{2}) \Psi^{*}_{-}(-q^{-\three}z_{1})
   \Big) \nonumber\\
   &&
\end{eqnarray}
At the values $z_{1}/z_{2} = q^2 (-p^{*\half})^{-\ell'}$ where
$\ell' \in \ZZ_{\ge 0}$, $\ell'$ odd, the function $1/\alpha^{*+}$ of the
r.h.s. of eq. (\ref{eq:splitpsiA}) has a zero while the function
$1/\alpha^{*+}$ of the l.h.s. has neither poles nor zeroes. 
at the values $z_{1}/z_{2} = q^2 (-p^{*\half})^{-\ell'}$ where
$\ell' \in \ZZ_{<0}$, $\ell'$ odd, the function $1/\alpha^{*+}$ of the
l.h.s. of eq. (\ref{eq:splitpsiA}) has a pole while the function
$1/\alpha^{*+}$ of the r.h.s. is regular.
It follows that the operator $\Psi^{*}_{+}(q^{-\half}z_{2})
\Psi^{*}_{+}(-q^{-\three}z_{1}) + \Psi^{*}_{-}(q^{-\half}z_{2})
\Psi^{*}_{-}(-q^{-\three}z_{1})$ must have a pole at $z_{1}/z_{2} =
q^2 (-p^{*\half})^{-\ell'}$, $\ell'+1 \in 2\ZZ$, 
while $\Psi^{*}_{+}(-q^{-\three}z_{1})
\Psi^{*}_{+}(q^{-\half}z_{2}) + \Psi^{*}_{-}(-q^{-\three}z_{1})
\Psi^{*}_{-}(q^{-\half}z_{2})$ is interpreted as a residue operator. \\
Finally, one concludes that the operator $T_{\ell\ell'}(z)$ has to be
interpreted as a product of residue operators from expression
(\ref{eq:T11vertex}) when both conditions $z_{1}/z_{2} =
(-p^{\half})^{-\ell}$ and $z_{1}/z_{2} = q^2 (-p^{*\half})^{-\ell'}$, with
$\ell$ and $\ell'$ odd, are fulfilled.

\medskip

One deals with the cases $\ell+\ell'$ odd
 along the same lines. Therefore, one can state:
\begin{theorem}\label{thm41}
  At $c=1$, the surface conditions $(-p^{\half})^{-\ell} = q^2 \,
  (-p^{*\half})^{-\ell'}$ where $\ell, \ell' \in \ZZ$, correspond to
  simultaneous existence of zeroes in the coefficients of the r.h.s.
  (when $\ell$ and/or $\ell' \in \ZZ_{\ge 0}$) or poles in the
  coefficients of the l.h.s. (when $\ell$ and/or $\ell' \in \ZZ_{<0}$) of
  the Riemann-Hilbert splitting of both products of vertex operators $\Phi$ 
  on the one hand
  and $\Psi^{*}$ on the other hand. \\
  The operators $T_{\ell\ell'}(z)$ are then interpreted as residue
  operators of $\cT_{\bar\ell\bar\ell'}(z_{1},z_{2})$ (expressed in term of
  vertex operators) at the point $z_{1}=\gamma_{\ell\ell'} z_{2}$.
\end{theorem}

Interpretation of the operators $S_{\ell\ell'}(z)$ follows along the same
lines, introducing
\begin{equation}
  \label{eq:defS12}
  \cS_{\bar\ell\bar\ell'}(z_{1},z_{2}) \equiv \tr \Big( \sigma_{x}^{\ell}
  \; L(z_{1}) \; \sigma_{x}^{\ell'} \; L(z_{2})^{-1} \Big) \qquad
  \ell\,,\ \ell'\in\ZZ
\end{equation}
which is related to the operator $S_{\ell\ell'}(z)$ by
$S_{\ell\ell'}(z)=\cS_{\bar\ell\bar\ell'}(z,q^{2}\gamma_{-\ell-\ell'}^{-1}z)$.
Performing the same analysis, one deduces similar conclusions for the
$S_{\ell\ell'}$ type operators, summarized in:
\begin{theorem}\label{thm42}
  At $c=1$, the surface conditions $(-p^{\half})^{\ell} = q^2 \,
  (-p^{*\half})^{\ell'}$ where $\ell, \ell' \in \ZZ$, correspond to
  simultaneous existence of zeroes in the coefficients of the r.h.s. (when
  $\ell$ and/or $\ell' \in \ZZ_{<0}$) or the poles of the coefficients of
  the l.h.s. (when $\ell$ and/or $\ell' \in \ZZ_{\ge 0}$) of the
  Riemann-Hilbert splitting of both products of vertex operators $\Phi$ on
  the one hand and $\Psi^{*}$ on the other hand. \\
  The operators $S_{\ell\ell'}(z)$ are then interpreted as residue
  operators of $\cS_{\bar\ell\bar\ell'}(z_{1},z_{2})$ (expressed in term of
  vertex operators) at the point $z_{1}=q^{-2}\gamma_{-\ell-\ell'} z_{2}$.
\end{theorem}

\section{Generalization to ${\cA}_{q,p}(\widehat{gl}_{N})$\label{sect:glN}}

\subsection{The algebra ${\cA}_{q,p}(\widehat{gl}_{N})$}

We start with the Boltzmann weights matrix for $\ZZ_{N}$-vertex model 
 \cite{Bela,ChCh}:
\begin{equation}
  \label{eq:RZN}
  \cW(z,q,p) = z^{2/N-2} \frac{1}{\kappa(z^2)}
  \frac{\displaystyle\vartheta\car{\sfrac{1}{2}}{\sfrac{1}{2}}(\mu,\tau)}
  {\displaystyle\vartheta\car{\sfrac{1}{2}}{\sfrac{1}{2}}(\xi+\mu,\tau)}
  \,\, \sum_{(\alpha_1,\alpha_2)\in\ZZ_N\times\ZZ_N}
  W_{(\alpha_1,\alpha_2)}(\xi,\mu,\tau) \,\, I_{(\alpha_1,\alpha_2)}
  \otimes I_{(\alpha_1,\alpha_2)}^{-1} 
\end{equation}
where the variables $z,q,p$ are related to the variables 
$\xi,\mu,\tau$ by
\begin{equation}
  z=e^{i\pi\xi} \,,\qquad q=e^{i\pi\mu} \,,\qquad p=e^{2i\pi\tau} 
\end{equation}
The Jacobi theta functions with rational characteristics
$\vartheta\car{\gamma_1}{\gamma_2}(\xi,\tau)$ are defined in Appendix 
\ref{app:theta}.
\\
The normalization factor is chosen as follows:
\begin{equation}
  \frac{1}{\kappa(z^2)} = \frac{(q^{2N}z^{-2};p,q^{2N})_\infty \,
  (q^2z^2;p,q^{2N})_\infty \, (pz^{-2};p,q^{2N})_\infty \,
  (pq^{2N-2}z^2;p,q^{2N})_\infty} {(q^{2N}z^2;p,q^{2N})_\infty \,
  (q^2z^{-2};p,q^{2N})_\infty \, (pz^2;p,q^{2N})_\infty \,
  (pq^{2N-2}z^{-2};p,q^{2N})_\infty} 
\end{equation}
The functions $W_{(\alpha_1,\alpha_2)}$ are given by
\begin{equation}
  W_{(\alpha_1,\alpha_2)}(\xi,\mu,\tau) =
  \frac{\displaystyle\vartheta\car{\sfrac{1}{2}+\alpha_1/N}
  {\sfrac{1}{2}+\alpha_2/N}(\xi+\mu/N,\tau)} {\displaystyle
  N\vartheta\car{\sfrac{1}{2}+\alpha_1/N}
  {\sfrac{1}{2}+\alpha_2/N}(\mu/N,\tau)} 
\end{equation}
The matrices $I_{(\alpha_1,\alpha_2)}$ are defined as follows:
\begin{equation}
  I_{(\alpha_1,\alpha_2)} = g^{\half} (g^{\alpha_2} \, h^{\alpha_1})
  g^{-\half} 
\end{equation}
where the $N \times N$ matrices $g$ and $h$ are given by $g_{ij} =
\omega^i\delta_{ij}$ and $h_{ij} = \delta_{i+1,j}$, the addition of 
indices being understood modulo $N$. They satisfy $hg = \omega gh$. \\
The matrix (\ref{eq:RZN}) is $\ZZ_N$-symmetric, that is
$\cW_{c+s\,,\,d+s}^{a+s\,,\,b+s} = \cW_{c\,,\,d}^{a\,,\,b}$ for any indices
$a,b,c,d,s \in \ZZ_N$ (the addition of indices being understood modulo $N$)
and the non-vanishing elements of the $\cW$ matrix are of the type
$\cW_{c\,,\,a+b}^{a\,,\,c+b}$.

\medskip

To define the elliptic quantum algebra
${\cA}_{q,p}(\widehat{gl}_{N})$, we introduce the following matrix, 
which differs from (\ref{eq:RZN}) by a
suitable normalization factor:
\begin{equation}
  R_{12}(z) \equiv R_{12}(z,q,p) = \tau_N(q^{\frac 12}z^{-1}) \,
  \cW_{12}(z,q,p)
  \label{eq:RtildeN}
\end{equation}
where the function $\tau_N(z)$ is defined by
\begin{equation}
  \tau_N(z) = z^{\frac{2}{N}-2} \,
  \frac{\Theta_{q^{2N}}(qz^2)}{\Theta_{q^{2N}}(qz^{-2})} 
  \label{eq:tauN}
\end{equation}
The function $\tau_N(z)$ is $q^{N}$-periodic and satisfies
$\tau_N(z) \, \tau_N(z^{-1}) = 1$. \\
The matrix $R_{12}$ is crossing-unitary \cite{RT86,ENSLAPP670}:
\begin{equation}
  \label{eq:crossunitN}
  \Big( R_{12}(z)^{t_2} \Big)^{-1} = \Big(
  R_{12}(q^Nz)^{-1} \Big)^{t_2}
\end{equation}
and obeys a quasi-periodicity property \cite{ENSLAPP670}:
\begin{eqnarray}
  R_{12}(-p^{\half}z) &=& (g^{\half} h g^{\half} \otimes 1)^{-1} \,
  R_{21}(z^{-1})^{-1} \, (g^{\half} h g^{\half} \otimes 1) \nonumber\\
  &=&
  \widetilde{\tau}_{N}(z)^{-1} \, (g^{\half} h g^{\half} \otimes 1) \,
  R_{12}(z) \, (g^{\half} h g^{\half} \otimes 1)
  \label{eq:quasiperN}
\end{eqnarray}
where 
\begin{equation}
    \label{eq:tautildeN}
    \widetilde{\tau}(z) = \tau(q^{\half}z^{-1}) \,
    \tau(q^{\half}z) = q^{2/N-2} \; \frac{\Theta_{q^{2N}}(q^2z^2)
    \; \Theta_{q^{2N}}(q^2z^{-2})}{\Theta_{q^{2N}}(z^2) \;
    \Theta_{q^{2N}}(z^{-2})}
\end{equation}
Hence, we get
\begin{equation}
  \label{eq:relFN}
  R_{12}\left((-p^{\half})^{-\ell}z \right) = F(\ell,z)
  \, ((g^{\half} h g^{\half})^\ell \otimes 1) \, R_{12}(z)
  \, ((g^{\half} h g^{\half})^\ell \otimes 1)
\end{equation}
where, for $\ell \geq 0$,
\begin{eqnarray}
  \label{eq:FllN}
  F(\ell,z) &=& \prod_{s=1}^\ell
  \widetilde{\tau}_{N}((-p^\half)^{-s}z)
\\
\label{eq:inverN}
  F(-\ell,z) &=& \frac{1}{F(\ell,(-p^{\half})^{\ell}z)} 
\end{eqnarray}

\medskip

We now define the elliptic quantum algebra ${\cA}_{q,p}(\widehat{gl}_{N})$
\cite{FIJKMY,JKOS} as an algebra of operators $L_{ij}(z)$ $\equiv
\sum_{n\in\ZZ} L_{ij}(n) \, z^n$ where $i,j \in \ZZ/N\ZZ$, encapsulated
into a $N \times N$ matrix
\begin{equation}
  L(z) = \left(
  \begin{array}{ccc} 
    L_{11}(z) & \cdots & L_{1N}(z) \\ 
    \vdots && \vdots \\ 
    L_{N1}(z) & \cdots & L_{NN}(z) \\ 
  \end{array}
  \right) 
\end{equation}
One defines ${\cal A}_{q,p}(\widehat{gl}(N)_c)$ by imposing the following
constraints on the $L_{ij}(z)$ (with the matrix $R_{12}$ given by eq.
(\ref{eq:RtildeN})):
\begin{equation}
  R_{12}(z_{1}/z_{2}) \, L_1(z_{1}) \, L_2(z_{2}) = L_2(z_{2}) \,
  L_1(z_{1}) \, R_{12}^{*}(z_{1}/z_{2})
\end{equation}
where $L_1(z) \equiv L(z) \otimes \II$, $L_2(z) \equiv \II \otimes L(z)$
and $R^{*}_{12}$ is defined by $R^{*}_{12}(z,q,p) \equiv
R_{12}(z,q,p^*=pq^{-2c})$. \\
The matrix $R^{*}_{12}$ obeys also the properties of crossing-unitarity
(\ref{eq:crossunitN}) and quasi-periodicity (\ref{eq:quasiperN}), this last
one being understood with the modified elliptic nome $p^*$. \\
The $q$-determinant $q$-$\det L(z)$ given by
\begin{equation}
  q\mbox{-}\det L(z) \equiv \sum_{\sigma\in{\mathfrak S}_N}
  \eps(\sigma) \prod_{i=1}^N L_{i,\sigma(i)}(z q^{i-N-1})
\end{equation}
($\eps(\sigma)$ being the signature of the permutation $\sigma$) is in the
center of ${\cal A}_{q,p}(\widehat{gl}(N)_c)$. It can be set to the value
$q^{\frac c2}$ so as to get
\begin{equation}
  {\cal A}_{q,p}(\widehat{sl}(N)_c) = {\cal A}_{q,p}(\widehat{gl}(N)_c)/
  \langle q\mbox{-}\det L - q^{\frac c2} \rangle 
\end{equation}

\subsection{Exchange algebras}

\begin{proposition}\label{prop:LTN}
  We define the operators $T_{\ell\ell'}(z)$ and
  $S_{\ell\ell'}(z)$, $\ell,\ell' \in \ZZ$, by
  \begin{eqnarray}
    \label{eq:operTN}
    T_{\ell\ell'}(z) &\!\!\equiv\!\!& \tr \Big( L(\gamma z)^{-1} \;
    (g^{\half} h g^{\half})^{\ell} \; L(z) \; (g^{\half} h
    g^{\half})^{\ell'} \Big) \nonumber \\
    &\!\!=\!\!& \tr \Big( (g^{\half} h^{t} g^{\half})^{\ell} \; {(L(\gamma
    z)^{-1})}^{t} \; (g^{\half} h^{t} g^{\half})^{\ell'} \; {L(z)}^{t}
    \Big) \\
    \label{eq:operSN}
    S_{\ell\ell'}(z) &\!\!\equiv\!\!& \tr \Big( (g^{\half} h
    g^{\half})^{\ell} \; L(z) \; (g^{\half} h g^{\half})^{\ell'} \;
    {(L(\gamma^{-1}q^{-N}z)^{-1})} \Big) \nonumber \\
    &\!\!=\!\!& \tr \Big( {L(z)}^{t} \; (g^{\half} h^{t} g^{\half})^{\ell}
    \; {(L(\gamma^{-1}q^{-N}z)^{-1})}^{t} \; (g^{\half} h^{t}
    g^{\half})^{\ell'} \Big)
  \end{eqnarray}
  The operators $T_{\ell\ell'}(z)$ have exchange relations with
  the generators $L(z)$ of ${\cA}_{q,p}(\widehat{gl}_{N})$:
  \begin{equation}
    L(z_2) \; T_{\ell\ell'}(z_1) = f_{\ell\ell'}(z_1/z_2)
    \; T_{\ell\ell'}(z_1) \; L(z_2)
  \end{equation}
  if the conditions 
  \begin{equation}
\gamma = (-p^{\half})^{-\ell} = (-p^{*\half})^{-\ell'}\, q^{N}
\end{equation}
   are fulfilled. \\
  In the same way, operators $S_{\ell\ell'}(z)$ have exchange
  relations with the generators $L(z)$ of
  ${\cA}_{q,p}(\widehat{gl}_{N})$:
  \begin{equation}
    L(z_2) \; S_{\ell\ell'}(z_1) = f_{\ell\ell'}(z_1/z_2)
    \; S_{\ell\ell'}(z_1) \; L(z_2)
  \end{equation}
  if the conditions 
  \begin{equation}
  \gamma = (-p^{\half})^{\ell}  =
  (-p^{*\half})^{\ell'}\, q^{N}
\end{equation}
  are fulfilled. \\
  The exchange function $f_{\ell\ell'}(z)$ is
  given by
  \begin{equation}
    \label{eq:functstrucTLN}
    f_{\ell\ell'}(z) = \frac{F^*(\ell',z)}{F(\ell,z)}
  \end{equation}
  where $F(\ell,z)$ is now given by (\ref{eq:FllN})--(\ref{eq:inverN}) and
  $F^{*}(\ell,z)$ is obtained from $F(\ell,z)$ by $p \to p^{*}$.
\end{proposition}
\textbf{Proof:} the proof is completely algebraic and follows exactly the
same lines as the one of Proposition \ref{prop:LT}.
\finproof
Following the steps of section \ref{sect:DVA}, we introduce:
\begin{definition}
In the parameters space $(q,p,c)$, the surface $\cP^{(N)}_{\ell\ell'}$ is
defined by the relation 
\begin{equation}
(-p^{\half})^{-\ell} = (-p^{*\half})^{-\ell'}\,q^{N}\,,\mbox{ with }
p^{*}=p\,q^{-2c}
\label{eq:surfN}
\end{equation}
The operators $T_{\ell\ell'}$ and $S_{-\ell,-\ell'}$ defined on this
surface are given by relations (\ref{eq:operTN}) and (\ref{eq:operSN}) with
$\gamma\equiv\gamma_{\ell\ell'}=(-p^{\half})^{-\ell}$.
\end{definition}

\begin{proposition}\label{prop:TTN}
When $c$ is a rational number, several surfaces $\cP^{(N)}_{\ell\ell'}$ can 
be defined simultaneously. \\
  On the surface $\cP^{(N)}_{\ell\ell'}\cap\cP^{(N)}_{\lambda\lambda'}$, the
  operators $T_{\ell\ell'}$ and $T_{\lambda\lambda'}$ satisfy the following
  exchange algebra
  \begin{equation}
    \label{eq:exchangeTTN}
    T_{\ell\ell'}(z_1) \; T_{\lambda\lambda'}(z_2) =
    \mathbf{F}_{\ell\ell'}^{\lambda\lambda'}(z_1/z_2) \;
    T_{\lambda\lambda'}(z_2) \; T_{\ell\ell'}(z_1)
  \end{equation}
  On the surface 
  $\cP^{(N)}_{-\ell,-\ell'}\cap\cP^{(N)}_{-\lambda,-\lambda'}$, the
  operators $S_{\ell\ell'}$ and $S_{\lambda\lambda'}$ satisfy the following
  exchange algebra 
  \begin{equation}
    \label{eq:exchangeSSN}
    S_{\ell\ell'}(z_1) \; S_{\lambda\lambda'}(z_2) =
    \mathbf{F}_{\ell\ell'}^{\lambda\lambda'}(z_1/z_2) \;
    S_{\lambda\lambda'}(z_2) \; S_{\ell\ell'}(z_1)
  \end{equation}
  On the surface $\cP^{(N)}_{\ell\ell'}\cap\cP^{(N)}_{-\lambda,-\lambda'}$, 
  the
  operators $T_{\ell\ell'}$ and $S_{\lambda\lambda'}$ satisfy the following
  exchange algebra 
  \begin{equation}
    \label{eq:exchangeTSN}
    T_{\ell\ell'}(z_1) \; S_{\lambda\lambda'}(z_2) = 
    \mathbf{F}_{\ell\ell'}^{\lambda\lambda'}(z_1/z_2) \;
    S_{\lambda\lambda'}(z_2) \; T_{\ell\ell'}(z_1)
  \end{equation}
  The exchange function $\mathbf{F}_{\ell\ell'}^{\lambda\lambda'}(z)$
  is given by
  \begin{equation}
    \label{eq:functstrucTTN}
    \mathbf{F}_{\ell\ell'}^{\lambda\lambda'}(z) =
    \frac{F(\ell,z)}{F^*(\ell',z)} \;
    \frac{F^*(\ell',\gamma_{\lambda\lambda'}^{-1}z)}
    {F(\ell,\gamma_{\lambda\lambda'}^{-1}z)}
  \end{equation}
  where $F(\ell,z)$ is given by (\ref{eq:FllN})--(\ref{eq:inverN}).
\end{proposition}
\textbf{Proof:} the proof follows the same lines as the one of Proposition
\ref{prop:TT}.

\medskip

Similarly to the $sl(2)$ case, we need to introduce a suitable
Riemann--Hilbert splitting in order to be able to define ordering relations
among the component operators $ T_{\ell\ell'}(z)$. It can be consistently 
chosen as
\begin{equation}
  \varphi_{\ell\ell'}(z_2/z_1) \; T_{\ell\ell'}(z_1) \;
  T_{\ell\ell'}(z_2) = \varphi_{\ell\ell'}(z_1/z_2) \;
  T_{\ell\ell'}(z_2) \; T_{\ell\ell'}(z_1)
\end{equation}
where
\begin{eqnarray}
  \varphi_{\ell\ell'}(x) &\!\!=\!\!& \frac{1}{(1-x^2)^{2|\ell|-2|\ell'|}} \;
  \prod_{s=1}^{|\ell|-1} \frac{(q^2 p^{-s} x^2;q^{2N})_\infty \, (q^{2N-2}
  p^{-s} x^2;q^{2N})_\infty} {(q^2 p^{s} x^2;q^{2N})_\infty \, (q^{2N-2}
  p^{s} x^2;q^{2N})_\infty} \; \Bigg[ \frac{(p^{s} x^2;q^{2N})_\infty}
  {(q^{2N} p^{-s} x^2;q^{2N})_\infty} \Bigg]^2 \; \nonumber \\
  & \!\!\times\!\!& \prod_{s=1}^{|\ell'|-1} \frac{(q^2 {p^*}^{s}
  x^2;q^{2N})_\infty \, (q^{2N-2} {p^*}^s x^2;q^{2N})_\infty} {(q^{2}
  {p^*}^{-s} x^2;q^{2N})_\infty \, (q^{2N-2} {p^*}^{-s} x^2;q^{2N})_\infty}
  \; \Bigg[ \frac{(q^{2N} {p^*}^{-s} x^2;q^{2N})_\infty} {({p^*}^{s}
  x^2;q^{2N})_\infty} \Bigg]^2
\end{eqnarray}
This allows one to properly define the algebra generated by
$T_{\ell\ell'}(z)$. The next step should now be to interpret possible
vertex operator representation of ${\cA}_{q,p}(\widehat{sl}_{N})$, so as to
give a characterization of the surfaces (\ref{eq:surfN}) on the same lines
as the one explicited in Theorems \ref{thm41} and \ref{thm42} for $sl(2)$.
To our knowledge, such a representation is not yet available.

\section{Conclusion}

We have established a connection between the vertex operator
representations available at $c=1$ and the so-called ``surface conditions''
characterizing the existence of subalgebras of $q$-$W_N$ type in
${\cA}_{q,p}(\widehat{sl}_{N})$, as defined in \cite{ENSLAPP649,ENSLAPP670,LAPTH690}.
We are now able to define new directions
of investigations which will either make use of, or further extend this
connection. It will undoubtedly lead to a better understanding of the $q$-$W_N$ algebra
structures which we have constructed.

Regarding the first option, it must be indicated here that the VO
construction used in \cite{Junichi04} leads to the identification of
very specific, so-called ``fused'' operators located at some precise
singular points of the VO exchange algebra. However, a deep
understanding of the connection between these non-singular fusion
locations and our own construction, is at present lacking.

The singular fused operators, by contrast, now allow us to get explicit
representations of our generators of $q$-$W_N$-type algebras, thereby
enabling us to better understand subtle properties of these algebras.
Noticeable amongst them is the existence of consistent central extensions.
This delicate question was touched upon in \cite{LAPTH680}, where central
extensions were built in a formal way by explicitly solving the coboundary
equations (Jacobi identity) for given $q$-Virasoro exchange functions. The
relevance, and explicit construction of these formal central extensions, is
however disputable: indeed the entanglement between the requirements of
normal ordering of $q$-$W_N$ generators (such as discussed in Section 3.3)
leading to Riemann-Hilbert splitting of the structure functions, and the
resolution of the consistency conditions for the central extensions
(depending also on these structure functions), is a delicate issue.
Hopefully the study of short-distance expansions of $q$-$W_N$ generators
using such explicit examples as are now available in \cite{Junichi04},
may shed light on this problem.

Regarding the second option, it is obvious that extensions of this
connection entail the comparison of our commensurability conditions with VO
representations both for $c \ne 1$ in the case of
${\cA}_{q,p}(\widehat{sl}_{2})$ and more generically for
${\cA}_{q,p}(\widehat{sl}_{N})$, $N \ge 3$. 

\appendix

\section{Jacobi theta functions with rational characteristics\label{app:theta}}

Let $\HH = \{ z\in\CC \,\vert\, \mbox{Im} z > 0 \}$ be the upper
half-plane and $\Lambda_\tau = \{ \lambda_1\,\tau + \lambda_2 \,\vert\,
\lambda_1,\lambda_2 \in \ZZ \,, \tau \in \HH \}$ the lattice with basis
$(1,\tau)$ in the complex plane. One sets $\omega = e^{2i\pi/N}$.

\medskip

One defines the Jacobi theta functions with rational characteristics
$\gamma=(\gamma_1,\gamma_2) \in \sfrac{1}{N} \ZZ \times \sfrac{1}{N} \ZZ$
by:
\begin{equation}
  \label{eq:A1}
  \vartheta\car{\gamma_1}{\gamma_2}(\xi,\tau) = \sum_{m \in \ZZ}
  \exp\Big(i\pi(m+\gamma_1)^2\tau + 2i\pi(m+\gamma_1)(\xi+\gamma_2) \Big)
\end{equation}
The functions $\displaystyle \vartheta\car{\gamma_1}{\gamma_2}(\xi,\tau)$
satisfy the following shift properties:
\begin{eqnarray}
  \label{eq:A2}
 \!\!\!\! &\!\!\!\!&\! \vartheta\car{\gamma_1+\lambda_1}{\gamma_2+\lambda_2}(\xi,\tau) =
  \exp\big(2i\pi \gamma_1\lambda_2\big) \,\,
  \vartheta\car{\gamma_1}{\gamma_2}(\xi,\tau) \\
 \!\!\!\!  &\!\!\!\!&\! \vartheta\car{\gamma_1}{\gamma_2}(\xi+\lambda_1\tau+\lambda_2,\tau) =
  \exp\big(-i\pi\lambda_1^2\tau-2i\pi\lambda_1\xi\big) \,
  \exp\big(2i\pi(\gamma_1\lambda_2 - \gamma_2\lambda_1)\big) \,
  \vartheta\car{\gamma_1}{\gamma_2}(\xi,\tau) \ \nonumber\\
  && \label{eq:A3}
\end{eqnarray}
where $\gamma=(\gamma_1,\gamma_2) \in \sfrac{1}{N}\ZZ \times
\sfrac{1}{N}\ZZ$ and $\lambda=(\lambda_1,\lambda_2) \in \ZZ \times \ZZ$. \\
Moreover, for arbitrary $\lambda=(\lambda_1,\lambda_2)$ (not necessarily
integers), one has the following shift exchange:
\begin{equation}
  \label{eq:A4}
  \vartheta\car{\gamma_1}{\gamma_2}(\xi+\lambda_1\tau+\lambda_2,\tau) =
  \exp\big(-i\pi\lambda_1^2\tau-2i\pi\lambda_1(\xi+\gamma_2+\lambda_2)\big)
  \, \vartheta\car{\gamma_1+\lambda_1}{\gamma_2+\lambda_2}(\xi,\tau) 
\end{equation}
The Jacobi theta functions
$\displaystyle\vartheta\car{\gamma_1}{\gamma_2}(\xi,\tau)$ with rational
characteristics $(\gamma_1,\gamma_2)$ can be expressed in terms of the
usual theta function $\Theta_p(z) = (z;p)_\infty \,
(pz^{-1};p)_\infty \, (p;p)_\infty$ as (with $p = e^{2i\pi\tau}$ and
$z = e^{i\pi \xi}$):
\begin{equation}
  \label{eq:A5}
  \vartheta\car{\gamma_1}{\gamma_2}(\xi,\tau) = e^{2i\pi\gamma_1\gamma_2}
  \, p^{\frac{1}{2}\gamma_1^2} \, z^{2\gamma_1} \,
  \Theta_{p}(-e^{2i\pi\gamma_2} p^{\gamma_1+\frac{1}{2}} z^2) 
\end{equation}

%%%%%%%%%%%%%%%%%%%%%%%%%%%%%%%%%%%%%%%%%%%%%%%%%%%%%%%%%%%%%%%%%%%%%%%%%%%

\end{document}